\documentclass[a4paper,12pt]{amsart}

\usepackage{float}
\usepackage{euscript,eufrak,verbatim}
\usepackage{graphicx}
\usepackage[usenames]{color}
\usepackage[colorlinks,linkcolor=red,anchorcolor=blue,citecolor=blue]{hyperref}
\usepackage{amsmath}
\usepackage{amsthm}
\usepackage[all]{xy}
\usepackage{graphicx} 
\usepackage{amssymb}

\usepackage{mathrsfs}
\usepackage{amscd}

\makeatletter
\newcommand{\svdots}{%
  \vbox{\fontsize{\sf@size}{\sf@size pt}\linespread{0.3}\selectfont
    \kern0.2\baselineskip
    \hbox{.}\hbox{.}\hbox{.}%
    \kern0.1\baselineskip
  }%
}
\makeatother

\theoremstyle{plain}
\newtheorem{main theorem}{Main Theorem}
\newtheorem{theorem}{Theorem}[section]
\newtheorem{lemma}[theorem]{Lemma}

\newtheorem{corollary}[theorem]{Corollary}
\newtheorem{proposition}[theorem]{Proposition}

\newtheorem{lemma-definition}[theorem]{Lemma-Definition}
\theoremstyle{definition}

\newtheorem{remark}[theorem]{Remark}

\makeatletter % `@' now normal "letter"
\@addtoreset{equation}{section}
\numberwithin{equation}{section}

\newcommand{\norm}[1]{\left\lVert#1\right\rVert}

\newcommand{\diam}{\mathrm{Diam}}

\newcommand{\umdimm}{\overline{\mathrm{mdim}}_{\mathrm{M}}}

\newcommand{\urdim}{\overline{\mathrm{rdim}}}
\newcommand{\lrdim}{\underline{\mathrm{rdim}}}
\newcommand{\rdim}{\mathrm{rdim}}

\title[Rate distortion dimension and ergodic decomposition]
{Rate distortion dimension and ergodic decomposition for $\mathbb{R}^d$-actions}

\author{Masaki Tsukamoto}

\address%[authorlabel1]
{Department of Mathematics, Kyoto University, Kitashirakawa Oiwake-cho, Sakyo-ku, Kyoto 606-8502, Japan}

\email{tsukamoto@math.kyoto-u.ac.jp}

\begin{document}



\thanks{The author was supported by JSPS KAKENHI JP21K03227.}

\maketitle

\begin{abstract}
Rate distortion dimension describes the theoretical limit of lossy data compression methods 
as the distortion bound goes to zero.
It was originally introduced in the context of information theory, and recently it was discovered that it 
has an intimate connection to Gromov’s theory of mean dimension of dynamical systems.
This paper studies the behavior of rate distortion dimension of $\mathbb{R}^d$-actions under ergodic decomposition.
Our main theorems provide natural convexity and concavity of upper and lower rate distortion dimensions under 
convex combination of invariant probability measures.
We also present examples which clarify the validity and limitations of the theorems.
\end{abstract}

\section{Main results}  \label{section: main results}

The purpose of this paper is to investigate the behavior of rate distortion dimension of $\mathbb{R}^d$-actions under ergodic decomposition.
Gromov \cite{Gromov} initiated the study of group actions on nonlinear infinite dimensional spaces arising from geometric analysis.
Typically such group actions have infinite topological entropy. 
Therefore standard entropy theory provides no useful information for them.
Gromov introduced a new topological invariant of group actions called \textit{mean dimension} for studying such objects.
Mean dimension is the \textit{averaged} number of parameters for describing orbits of the given group actions. 
This provides a useful information for infinite dimensional and infinite entropy group actions.
Mean dimension theory has been studied for more than 20 years and several applications have been discovered 
\cite{Lindenstrauss--Weiss, Lindenstrauss, Gutman--Lindenstrauss--Tsukamoto, Gutman--Tsukamoto_embedding, Gutman--Qiao--Tsukamoto}.

For developing broader applications of mean dimension, Lindenstrauss and the author 
\cite{Lindenstrauss--Tsukamoto_IEEE, Lindenstrauss--Tsukamoto_double} introduced a variational principle for mean dimension.
(So far, several authors have introduced their own approaches to variational principles in mean dimension theory; 
see e.g. \cite{Gutman--Spiewak, Shi22}.)
The theory of \cite{Lindenstrauss--Tsukamoto_IEEE, Lindenstrauss--Tsukamoto_double} connects mean dimension to 
a quantity called \textit{rate distortion dimension} \cite{Kawabata--Dembo}.
This was first introduced in the context of rate distortion theory, which is an important branch of information theory 
describing the theoretical limit of lossy  data compression schemes.
A familiar example of lossy data compression is JPEG algorithm for two-dimensional image compression.
Given an image, JPEG expands it by a wavelet basis and discards small coefficients.
Then the amount of information is significantly reduced.
Rate distortion theory provides a theoretical limit of such data compression scheme.

Recently the author \cite{Tsukamoto_random_Brody_curves} applied the variational principle of mean dimension theory 
to the study of entire holomorphic curves.
He found that they satisfy an inequality analogous to the \textit{Ruelle inequality} of smooth ergodic theory.
This opens an unexpected new interaction between geometric analysis and hyperbolic dynamics.
Such a research direction is probably close to the original spirit of Gromov.

For further developing applications of the variational principle, we need to know several basic properties of rate distortion dimension.
In particular it is fundamental to investigate its behavior under \textit{ergodic decomposition}.
Every invariant probability measure can be decomposed into convex combination of ergodic measures, 
and it is always important to understand the behavior of dynamical quantities under ergodic decomposition.
The purpose of the present paper is to study it for rate distortion dimension of $\mathbb{R}^d$-actions.
Maybe some readers wonder why we concentrate on the group $\mathbb{R}^d$. 
There are many other groups. Is there any good reason to study only $\mathbb{R}^d$?
The reason is that this is the most basic case for geometric applications.
For example, the study of entire holomorphic curves in the paper \cite{Tsukamoto_random_Brody_curves} 
provides an example of $\mathbb{R}^2$-actions.
We also notice that our results hold for $\mathbb{Z}^d$-actions as well 
(indeed, the case of $\mathbb{Z}^d$ is a bit simpler than that of $\mathbb{R}^d$), 
but a generalization to noncommutative groups is out of scope of this paper\footnote{This is mainly because the author does not have 
good examples motivating a generalization to noncommutative groups. Interested readers may pursue such a direction.}.

Throughout of the paper we assume that $d$ is a positive integer and that $\mathbb{R}^d$ is endowed with the Euclidean topology and 
standard additive group structure.
Let $(\mathcal{X}, \mathbf{d})$ be a compact metric space. 
Let $T\colon \mathbb{R}^d\times \mathcal{X} \to \mathcal{X}$ be a continuous action of $\mathbb{R}^d$ on $\mathcal{X}$.
A Borel probability measure $\mu$ on $\mathcal{X}$ is said to be \textbf{$T$-invariant} if 
$\mu\left(T^{-u}A\right) = \mu(A)$ for all $u\in \mathbb{R}^d$ and all Borel subsets $A\subset \mathcal{X}$.
We denote by $\mathscr{M}^T(\mathcal{X})$ the set of all $T$-invariant Borel probability measures on $\mathcal{X}$.
This space is endowed with the weak$^*$ topology.
A measure $\mu\in \mathscr{M}^T(\mathcal{X})$ is said to be \textbf{ergodic} if there is no $T$-invariant\footnote{A set $A\subset \mathcal{X}$ is 
said to be $T$-invariant if $T^{-u}A = A$ for all $u\in \mathbb{R}^d$.} Borel subset $A\subset \mathcal{X}$
with $0<\mu(A)<1$. 
We define $\mathscr{M}^T_{\mathrm{erg}}(\mathcal{X})\subset \mathscr{M}^T(\mathcal{X})$ as the set of all ergodic measures. 

Let $\mu\in \mathscr{M}^T(\mathcal{X})$. 
We randomly choose a point $x\in \mathcal{X}$ according to $\mu$ and consider its orbit 
$\{T^u x\}_{u\in \mathbb{R}^d}$.
For $\varepsilon>0$, we define the \textit{rate distortion function} $R(\mathbf{d}, \mu, \varepsilon)$
as the number of bits per unit volume of $\mathbb{R}^d$ for describing the orbit $\{T^u x\}_{u\in \mathbb{R}^d}$
within averaged distortion (w.r.t. $\mathbf{d}$) bounded by $\varepsilon$.
We will give the precise definition of $R(\mathbf{d}, \mu, \varepsilon)$ in \S \ref{subsection: rate distortion function}.
We define the \textbf{upper and lower rate distortion dimensions} by 
\[  \urdim(\mathcal{X}, T, \mathbf{d}, \mu) = \limsup_{\varepsilon\to 0} \frac{R(\mathbf{d}, \mu, \varepsilon)}{\log(1/\varepsilon)},\quad 
    \lrdim(\mathcal{X}, T, \mathbf{d}, \mu) = \liminf_{\varepsilon\to 0} \frac{R(\mathbf{d}, \mu, \varepsilon)}{\log(1/\varepsilon)}. \]
These are originally introduced by Kawabata--Dembo \cite{Kawabata--Dembo} 
in the context of information theory\footnote{The motivation of Kawabata and Dembo is to study 
how much information is carried by a signal that takes values in fractal sets.}.
When the upper and lower limits coincide, we denote the common value by 
$\rdim\left(\mathcal{X}, T, \mathbf{d}, \mu\right)$.    
The subject of this paper is to study the behavior of these quantities under ergodic decomposition 
and, more generally, convex combination of 
invariant probability measures.

For explaining our main results,
we need to introduce one more quantity called \textit{upper metric mean dimension}.
Let $(\mathcal{X}, \mathbf{d})$ be a compact metric space.
For $\varepsilon>0$ we define the \textbf{$\varepsilon$-covering number} 
$\#\left(\mathcal{X}, \mathbf{d}, \varepsilon\right)$ as the minimum integer $n>0$ for which there 
exists an open covering $\mathcal{X} = U_1\cup U_2\cup \dots \cup U_n$ with $\diam\left(U_i, \mathbf{d}\right) < \varepsilon$
for all $1\leq i \leq n$. Here $\diam \left(U_i, \mathbf{d}\right) := \sup_{x, y\in U_i} \mathbf{d}(x, y)$.
Let $T\colon \mathbb{R}^d\times \mathcal{X}\to \mathcal{X}$ be a continuous action of $\mathbb{R}^d$ on $\mathcal{X}$.
For $L>0$ we define a new metric $\mathbf{d}_L$ on $\mathcal{X}$ by 
\[  \mathbf{d}_L(x, y) = \sup_{u\in [0, L)^d} \mathbf{d}(T^u x, T^u y). \]
We define the \textbf{entropy at the scale $\varepsilon$} by 
\begin{equation} \label{eq: entropy at the scale epsilon}
   S\left(\mathcal{X}, T, \mathbf{d}, \varepsilon\right) 
  = \lim_{L\to \infty} \frac{\log \#\left(\mathcal{X}, \mathbf{d}_L, \varepsilon\right)}{L^d}. 
\end{equation}  
This limit always exists (and its value is finite).
Finally we define the \textbf{upper metric mean dimension} by 
\[  \umdimm \left(\mathcal{X}, T, \mathbf{d}\right) = \limsup_{\varepsilon\to 0} 
   \frac{S \left(\mathcal{X}, T, \mathbf{d}, \varepsilon\right)}{\log (1/\varepsilon)}.   \]
This was originally introduced by Lindenstrauss--Weiss \cite{Lindenstrauss--Weiss}.
It is easy to see that the upper metric mean dimension is an upper bound on rate distortion dimension:
\[  \lrdim\left(\mathcal{X}, T, \mathbf{d}, \mu\right) \leq  \urdim\left(\mathcal{X}, T, \mathbf{d}, \mu\right)
     \leq \umdimm\left(\mathcal{X}, T, \mathbf{d}\right)  \quad (\forall \mu\in \mathscr{M}^T(\mathcal{X})). \]

Now we can state our main results.

\begin{theorem}  \label{theorem: upper rate distortion dimension and convex combination}
Let $T\colon \mathbb{R}^d\times \mathcal{X}\to \mathcal{X}$ be a continuous action of $\mathbb{R}^d$ on a compact metric space 
$(\mathcal{X}, \mathbf{d})$.
Suppose the upper metric mean dimension $\umdimm\left(\mathcal{X}, T, \mathbf{d}\right)$ is finite.
Let $\lambda$ be a Borel probability measure on $\mathscr{M}^T\left(\mathcal{X}\right)$ and define 
$\mu \in \mathscr{M}^T\left(\mathcal{X}\right)$ by 
\[  \mu = \int_{\mathscr{M}^T\left(\mathcal{X}\right)} \nu \, d\lambda(\nu). \]
Then we have 
\begin{equation} \label{eq: convexity of upper rate distortion dimension}
    \urdim\left(\mathcal{X}, T, \mathbf{d}, \mu\right) \leq \int_{\mathscr{M}^T\left(\mathcal{X}\right)} 
     \urdim\left(\mathcal{X}, T, \mathbf{d}, \nu\right) \, d\lambda(\nu). 
\end{equation}     
\end{theorem}

\begin{theorem}  \label{theorem: lower rate distortion dimension and convex combination}
Let $T\colon \mathbb{R}^d\times \mathcal{X}\to \mathcal{X}$ be a continuous action of $\mathbb{R}^d$ on a compact metric space 
$(\mathcal{X}, \mathbf{d})$.
Let $\lambda$ be a Borel probability measure on $\mathscr{M}^T\left(\mathcal{X}\right)$ and define 
$\mu \in \mathscr{M}^T\left(\mathcal{X}\right)$ by 
\[  \mu = \int_{\mathscr{M}^T\left(\mathcal{X}\right)} \nu \, d\lambda(\nu). \]
Then we have 
 \begin{equation} \label{eq: concavity of lower rate distortion dimension}
     \lrdim\left(\mathcal{X}, T, \mathbf{d}, \mu\right) \geq \int_{\mathscr{M}^T\left(\mathcal{X}\right)} 
        \lrdim\left(\mathcal{X}, T, \mathbf{d}, \nu\right) \, d\lambda(\nu). 
  \end{equation}      
\end{theorem}

We will see later that the inequalities (\ref{eq: convexity of upper rate distortion dimension}) and 
(\ref{eq: concavity of lower rate distortion dimension}) may be strict in general.
(See the discussion after Theorem \ref{theorem: irregular behavior of upper rate distortion dimension} below and 
\S \ref{section: another example}.)

Notice that the condition $\umdimm\left(\mathcal{X}, T, \mathbf{d}\right) < \infty$ is required in Theorem 
\ref{theorem: upper rate distortion dimension and convex combination} whereas it is not assumed in 
Theorem \ref{theorem: lower rate distortion dimension and convex combination}.
The condition $\umdimm\left(\mathcal{X}, T, \mathbf{d}\right) < \infty$ is not a severe restriction.
As far as the author knows, it is satisfied for all interesting geometric examples
which have been studied so far.
By combining Theorems \ref{theorem: upper rate distortion dimension and convex combination} and 
\ref{theorem: lower rate distortion dimension and convex combination}, we get a corollary:

\begin{corollary} \label{corollary: rate distortion dimension and convex combination}
Let $T\colon \mathbb{R}^d\times \mathcal{X}\to \mathcal{X}$ be a continuous action of $\mathbb{R}^d$ on a compact metric space 
$(\mathcal{X}, \mathbf{d})$ with $\umdimm\left(\mathcal{X}, T, \mathbf{d}\right) < \infty$.
Let $\lambda$ be a Borel probability measure on $\mathscr{M}^T\left(\mathcal{X}\right)$ and define 
$\mu \in \mathscr{M}^T\left(\mathcal{X}\right)$ by 
\[  \mu = \int_{\mathscr{M}^T\left(\mathcal{X}\right)} \nu \, d\lambda(\nu). \]
If $\rdim\left(\mathcal{X}, T, \mathbf{d}, \nu\right)$ exists for $\lambda$-almost every $\nu\in \mathscr{M}^T(\mathcal{X})$,
then $\rdim\left(\mathcal{X}, T, \mathbf{d}, \mu\right)$ also exists and satisfies 
\[ \rdim\left(\mathcal{X}, T, \mathbf{d}, \mu\right) = \int_{\mathscr{M}^T\left(\mathcal{X}\right)} 
        \rdim\left(\mathcal{X}, T, \mathbf{d}, \nu\right) \, d\lambda(\nu). \] 
\end{corollary}

The next result shows that we cannot remove the assumption $\umdimm\left(\mathcal{X}, T, \mathbf{d}\right) < \infty$
 in Theorem \ref{theorem: upper rate distortion dimension and convex combination}.

\begin{theorem} \label{theorem: irregular behavior of upper rate distortion dimension}
There exist a compact metric space $(\mathcal{X}, \mathbf{d})$ and a continuous action 
$T\colon \mathbb{R}^d \times \mathcal{X} \to \mathcal{X}$ satisfying the following two conditions.
   \begin{itemize}
     \item For any ergodic measure $\nu \in \mathscr{M}^T_{\mathrm{erg}}\left(\mathcal{X}\right)$ we have 
              $\rdim\left(\mathcal{X}, T, \mathbf{d}, \nu\right) = 0$.
     \item There exists $\mu\in \mathscr{M}^T\left(\mathcal{X}\right)$ with $\rdim\left(\mathcal{X}, T, \mathbf{d}, \mu\right) = \infty$.         
   \end{itemize}
\end{theorem}

Let $\mu$ be the measure given in the second condition of this theorem.
By the ergodic decomposition theorem, 
there exists a Borel probability measure $\lambda$ on $\mathscr{M}^T_{\mathrm{erg}}\left(\mathcal{X}\right)$
for which we have
\[  \mu = \int_{\mathscr{M}^T_{\mathrm{erg}}\left(\mathcal{X}\right)} \nu\, d\lambda(\nu). \]
Then it follows from the first condition of the theorem that
\[  \int_{\mathscr{M}^T_{\mathrm{erg}}\left(\mathcal{X}\right)} \urdim\left(\mathcal{X}, T, \mathbf{d}, \nu\right) d\lambda(\nu) = 0, \]
while we have $\rdim\left(\mathcal{X}, T, \mathbf{d}, \mu\right) = \infty$.
Therefore the inequality (\ref{eq: convexity of upper rate distortion dimension}) in 
Theorem \ref{theorem: upper rate distortion dimension and convex combination}
does not hold here. (The system $(\mathcal{X}, T,\mathbf{d})$ has infinite upper metric mean dimension.
So it does not conflict with Theorem \ref{theorem: upper rate distortion dimension and convex combination}.)
Hence Theorem \ref{theorem: irregular behavior of upper rate distortion dimension} 
shows that we cannot remove the assumption on the upper metric mean dimension in Theorem 
\ref{theorem: upper rate distortion dimension and convex combination}.
It also shows that 
the inequality (\ref{eq: concavity of lower rate distortion dimension}) 
in Theorem \ref{theorem: lower rate distortion dimension and convex combination}
can be a strict inequality in general.

The main ingredient of the proofs of Theorems \ref{theorem: upper rate distortion dimension and convex combination} and 
\ref{theorem: lower rate distortion dimension and convex combination} is the formula of rate distortion function 
under convex combination of measures
(Theorem \ref{theorem: rate distortion function and convex combination of measures} 
in \S \ref{subsection: rate distortion function and convex combination of measures}):
Let $T\colon \mathbb{R}^d\times \mathcal{X}\to \mathcal{X}$ be a continuous action on a compact metric space 
$(\mathcal{X}, \mathbf{d})$.
Let $\lambda$ be a Borel probability measure on $\mathscr{M}^T\left(\mathcal{X}\right)$ and set
$\mu := \int_{\mathscr{M}^T(\mathcal{X})} \nu\, d\lambda(\nu)$.
Then for any $\varepsilon>0$ 
\begin{equation} \label{eq: rate distortion function and convex combination}
   R(\mathbf{d}, \mu, \varepsilon)
   = \inf\left\{\int_{\mathscr{M}^T(\mathcal{X})} R(\mathbf{d},\nu, \varepsilon_\nu)d\lambda(\nu) \middle|\, 
\text{\parbox{7cm}{$\mathscr{M}^T(\mathcal{X})\ni \nu \mapsto \varepsilon_\nu \in (0, \infty)$ with  
$\int_{\mathscr{M}^T(\mathcal{X})} \varepsilon_\nu \, d\lambda(\nu) \leq \varepsilon$}}\right\}.
\end{equation}  
Here the infimum is taken over all measurable maps $\mathscr{M}^T(\mathcal{X})\ni \nu \mapsto \varepsilon_\nu \in (0, \infty)$
satisfying $\int_{\mathscr{M}^T(\mathcal{X})} \varepsilon_\nu \, d\lambda(\nu) \leq \varepsilon$.
Theorems \ref{theorem: upper rate distortion dimension and convex combination} and 
\ref{theorem: lower rate distortion dimension and convex combination} follow from this formula.

The equation (\ref{eq: rate distortion function and convex combination})
is well-known for $\mathbb{Z}$-actions in 
classical information theory literature \cite{Gray--Davisson, Kieffer75, Shields--Neuhoff--Davisson--Ledrappier,
LeonGarcia--Davisson--Neuhoff, Effros--Chou--Gray}.
The author cannot find a paper studying its generalization to $\mathbb{R}^d$-actions, 
and it is not very easy for pure mathematicians, including the author himself, to read 
the papers \cite{Gray--Davisson, Kieffer75, Shields--Neuhoff--Davisson--Ledrappier,
LeonGarcia--Davisson--Neuhoff, Effros--Chou--Gray}.
(However, it is instructive to read them because they provide broader and 
more coherent perspective about rate distortion theory 
for nonergodic sources\footnote{Probably the introduction of \cite{LeonGarcia--Davisson--Neuhoff}
is the best one for layman to understand basic ideas.}.)
So we provide the full proof of (\ref{eq: rate distortion function and convex combination}).
This is a purely technical task and there is no new idea in our argument.
Therefore this paper is primally a technical paper and do not contain any innovation.
However, Theorems \ref{theorem: upper rate distortion dimension and convex combination} and 
\ref{theorem: lower rate distortion dimension and convex combination} seem to be fundamental, 
and it is desirable to publish the detailed proofs.

\begin{remark}
The present paper concentrates on the study of $\mathbb{R}^d$-actions.
Probably $\mathbb{Z}^d$-actions are another important class of group actions in geometric applications.
The statements of Theorems \ref{theorem: upper rate distortion dimension and convex combination},
\ref{theorem: lower rate distortion dimension and convex combination}, \ref{theorem: irregular behavior of upper rate distortion dimension}
and Corollary \ref{corollary: rate distortion dimension and convex combination} hold for $\mathbb{Z}^d$-actions as well.
(We have not found a paper containing these statements even for $\mathbb{Z}$-actions.)
Since the case of $\mathbb{Z}^d$ is technically simpler than that of $\mathbb{R}^d$, 
we provide the full arguments for the case of $\mathbb{R}^d$-actions. 
\end{remark}

\section{Mutual information and rate distortion function}  

The purpose of this section is to prepare basic definitions and properties of 
mutual information and rate distortion function.
The readers can find more details in \cite[Section 2]{Tsukamoto_Variational_R^d}.
Throughout the present paper, we assume that the base of logarithm is two:
\[  \log x = \log_2 x. \]

\subsection{Measure theoretic details}  \label{subsection: measure theoretic details}

Here we prepare some facts on measure theory.
A pair $(\mathcal{X}, \mathcal{A})$ is called a \textbf{measurable space} if 
$\mathcal{X}$ is a set and $\mathcal{A}$ is its $\sigma$-algebra.
A triplet $(\mathcal{X}, \mathcal{A}, \mathbb{P})$ is called a \textbf{probability space} if 
$(\mathcal{X}, \mathcal{A})$ is a measurable space and $\mathbb{P}$ is a probability measure on it.

We always assume that a topological space is equipped with its \textbf{Borel $\sigma$-algebra}
(the smallest $\sigma$-algebra containing all open sets).
We also assume that a finite set is equipped with the discrete topology and discrete $\sigma$-algebra
(the set of all subsets).
A topological space $\mathcal{X}$ is called a \textbf{Polish space} 
if it admits a metric $\mathbf{d}$ compatible with its given topology for which
$(\mathcal{X}, \mathbf{d})$ is a complete separable metric space.
A measurable space $(\mathcal{X}, \mathcal{A})$ is called a \textbf{standard Borel space} if it is isomorphic 
(as a measurable space) to some Polish space equipped with its Borel $\sigma$-algebra. 
Readers can find basic information about standard Borel spaces in the book of \cite{Srivastava}.

Let $(\mathcal{X}, \mathcal{A})$ and $(\mathcal{Y}, \mathcal{B})$ be measurable spaces.
A map $\nu\colon \mathcal{X}\times \mathcal{B} \to [0,1]$ is called a 
\textbf{transition probability on $\mathcal{X}\times \mathcal{Y}$}
if it satisfies the following two conditions.
\begin{itemize}
   \item For every $x\in \mathcal{X}$, the map $\mathcal{B}\ni B \mapsto \nu(x, B)\in [0,1]$ is a probability measure on 
   $(\mathcal{Y}, \mathcal{B})$.
   \item  For every $B\in \mathcal{B}$, the map $\mathcal{X}\ni x \mapsto \nu(x, B) \in [0,1]$ is measurable.
\end{itemize} 
We often denote $\nu(x, B)$ by $\nu(B|x)$.
We define the \textbf{product $\sigma$-algebra} $\mathcal{A}\otimes \mathcal{B}$ as the smallest $\sigma$-algebra of $\mathcal{X}\times \mathcal{Y}$
containing all ``rectangles” $A\times B$ $(A\in \mathcal{A}, B\in \mathcal{B})$.
For any $E\in \mathcal{A}\otimes \mathcal{B}$ and $x\in \mathcal{X}$, 
the section $E_x := \{y\in \mathcal{Y}\mid (x, y)\in E\}$ belongs to $\mathcal{B}$.
Moreover, if $(\mathcal{Y}, \mathcal{B})$ is a standard Borel space, then 
the map $\mathcal{X}\ni x\mapsto \nu(E_x|x)\in [0,1]$ is measurable for any transition probability 
$\nu$ on $\mathcal{X}\times \mathcal{Y}$ \cite[Proposition 3.4.24]{Srivastava}.

Let $(\Omega, \mathcal{F}, \mathbb{P})$ be a probability space and $(\mathcal{X}, \mathcal{A})$ a measurable space.
For a measurable map $X\colon \Omega \to \mathcal{X}$, we denote the push-forward measure $X_*\mathbb{P}$ by 
$\mathrm{Law} X$ and call it \textbf{the law of $X$} or \textbf{the distribution of $X$}.
Here $X_*\mathbb{P}$ is a probability measure on $\mathcal{X}$ defined by $X_*\mathbb{P}(A) = \mathbb{P}\left(X^{-1}(A)\right)$.

The next theorem introduces the notion of regular conditional distribution.
For the details, see \cite[p.15 Theorem 3.3 and Corollary]{Ikeda--Watanabe} or \cite[p.182 Corollary 6.2]{Gray_probability}.

\begin{theorem}  \label{theorem: regular conditional probability}
Let $(\Omega, \mathcal{F}, \mathbb{P})$ be a probability space.
Let $(\mathcal{X}, \mathcal{A})$ and $(\mathcal{Y}, \mathcal{B})$ be standard Borel spaces, and let 
$X\colon \Omega \to \mathcal{X}$ and $Y\colon \Omega \to \mathcal{Y}$ be measurable maps.
Set $\mu = \mathrm{Law} X$.
Then there exists a transition probability $\nu$ on $\mathcal{X}\times \mathcal{Y}$ such that for any $E\in \mathcal{A}\otimes \mathcal{B}$
we have 
\[ \mathbb{P}\left((X, Y)\in E\right) = \int_{\mathcal{X}} \nu(E_x|x) \, d\mu(x). \]
Moreover, the transition probability $\nu$ is unique in the following sense: 
If another transition probability $\nu^\prime$ on 
$\mathcal{X}\times \mathcal{Y}$ satisfies the same property then there exists a $\mu$-null set $N\in \mathcal{A}$ such that 
$\nu(B|x) = \nu^\prime(B|x)$ for all $x\in \mathcal{X}\setminus N$ and $B\in \mathcal{B}$.
\end{theorem}

The transition probability $\nu$ introduced in this theorem is called 
\textbf{the regular conditional distribution of $Y$ given $X=x$}.
We often denote $\nu(B|x)$ by $\mathbb{P}\left(Y\in B|X=x\right)$.

Let $(\mathcal{X}, \mathcal{A}), (\mathcal{Y},\mathcal{B}), (\mathcal{Z},\mathcal{C})$ be standard Borel spaces, and
let $X, Y, Z$ be random variables (defined on some common probability space $(\Omega, \mathcal{F}, \mathbb{P})$) 
taking values in $\mathcal{X}, \mathcal{Y}, \mathcal{Z}$ respectively.
We say that $X$ and $Z$ are \textbf{conditionally independent given $Y$} if there exists a $Y_*\mathbb{P}$-null set 
$N\subset \mathcal{Y}$ such that 
\[ \mathbb{P}\left((X, Z)\in A\times C \middle| Y=y\right) = 
    \mathbb{P}(X\in A|Y=y) \cdot \mathbb{P}(Z\in C|Y=y) \]
for all $A\in \mathcal{A}, C \in \mathcal{C}$ and $y \in \mathcal{Y} \setminus N$.
Here the left-hand side is the regular conditional distribution of $(X, Z)$ given $Y=y$.
In information theory literature \cite[p.34]{Cover--Thomas}, 
one says that random variables $X, Y, Z$ \textbf{form a Markov chain in this order}
(denoted by $X\to Y \to Z$) if $X$ and $Z$ are conditionally independent given $Y$.
This notion will be important in the \textit{data-processing inequality}
(Lemma \ref{lemma: data processing inequality}) below.

\subsection{Mutual information}  \label{subsection: mutual information}
Here we prepare basic definitions and properties of mutual information.
Readers can find a nice introductory exposition in \cite{Cover--Thomas}.
The book of Gray \cite{Gray_entropy} provides mathematically sophisticated details of the theory.
 
Throughout this subsection, we fix a probability space $(\Omega, \mathcal{F}, \mathbb{P})$ and assume that
all random variables are defined on it.
For a random variable $X$ that takes values in a finite set $A$, we define its \textbf{Shannon entropy} by
\[  H(X) = -\sum_{a\in A} \mathbb{P}(X=a) \log \mathbb{P}(X=a). \]
Here we assume $0\log 0  =0$.

Let $X$ and $Y$ be random variables that take values in measurable space $(\mathcal{X}, \mathcal{A})$ and 
$(\mathcal{Y}, \mathcal{B})$ respectively.
We would like to define the \textbf{mutual information} $I(X;Y)$ which estimates the amount of information shared by 
$X$ and $Y$.
We need to consider the two cases:

\begin{enumerate}
    \item[(1)] Suppose $\mathcal{X}$ and $\mathcal{Y}$ are finite sets. Then we define 
    \[  I(X;Y) = H(X) + H(Y) - H(X, Y). \]
    Here $H(X, Y)$ is the Shannon entropy of the pair $(X, Y)\colon \Omega \to \mathcal{X}\times \mathcal{Y}$.
    We always have $I(X; Y) \leq \min\left(H(X), H(Y)\right)$. 
    If $f\colon \mathcal{X}\to \mathcal{X}^\prime$ and $g\colon \mathcal{Y}\to \mathcal{Y}^\prime$ are maps with 
    finite sets $\mathcal{X}^\prime$ and $\mathcal{Y}^\prime$, then we have 
    \begin{equation}  \label{eq: pre-data processing inequality}
       I\left(f(X); g(Y)\right) \leq I\left(X; Y\right). 
    \end{equation}   
    \item[(2)] In general, we consider finite measurable partitions $\alpha = \{A_1, \dots, A_m\}$ and 
    $\beta = \{B_1, \dots, B_n\}$ of $\mathcal{X}$ and $\mathcal{Y}$ respectively. 
    We define a random variable $\alpha\circ X\colon \Omega \to \{1, 2, \dots, m\}$ by $\alpha\circ X = i \Longleftrightarrow  X\in A_i$.
    We define $\beta \circ Y$ similarly.
    These random variables take only finitely many values. So we can consider their mutual information
    $I\left(\alpha\circ X; \beta\circ Y\right)$ by the formula
    \[ I(\alpha\circ X; \beta \circ Y) = H\left(\alpha\circ X\right) + H\left(\beta\circ Y\right) - H\left(\alpha\circ X, \beta\circ Y\right). \]
    We define the mutual information $I(X; Y)$ as the supremum of $I(\alpha\circ X; \beta\circ Y)$
    over all finite measurable partitions $\alpha$ and $\beta$ of $\mathcal{X}$ and $\mathcal{Y}$.
    It follows from the inequality (\ref{eq: pre-data processing inequality}) that this definition is compatible with the case (1). 
    Mutual information is symmetric and nonnegative: $I(X; Y) = I(Y; X) \geq 0$.
    If $X$ and $Y$ are independent then $I(X; Y) = 0$.
\end{enumerate}    

If $(\mathcal{X}, \mathcal{A})$ and $(\mathcal{Y}, \mathcal{B})$ are standard Borel spaces, 
then we can provide a slightly different 
description of $I(X; Y)$.
Let $\mu =\mathrm{Law} X$ be the distribution of $X$, and let 
$\nu(B|x) = \mathbb{P}(Y\in B|X=x)$ $(B\in \mathcal{B}, x\in \mathcal{X})$
be the regular conditional distribution of $Y$ given $X=x$.
The (joint) distribution of $(X, Y)$ is determined by $\mu$ and $\nu$.
In particular $I(X; Y)$ is also determined by them.
So we sometimes denote $I(X; Y)$ by $I(\mu, \nu)$.
The usefulness of this description is that $I(\mu, \nu)$ is a concave function in $\mu$ and a convex function in $\nu$
(see Proposition \ref{proposition: concavity and convexity of mutual information} below).

In some arguments we also need to use \textit{conditional mutual information}.
Let $X, Y, Z$ be random variables taking values in standard Borel spaces 
$(\mathcal{X}, \mathcal{A}), (\mathcal{Y}, \mathcal{B}), (\mathcal{Z},\mathcal{C})$ respectively.
We would like to define the conditional mutual information $I(X; Y|Z)$.
Let $\lambda = Z_*\mathbb{P}$ be the distribution of $Z$.
For each $z\in \mathcal{Z}$ we define a probability measure $\mu_z$ on $\mathcal{X}$ and 
a transition probability $\nu_z$ on $\mathcal{X}\times \mathcal{Y}$ by 
\begin{align*}
    \mu_z(A) & = \mathbb{P}(X\in A|Z=z) \quad (A\in \mathcal{A}), \\
    \nu_z(B|x) & = \mathbb{P}\left(Y\in B|(X, Z) = (x, z)\right) \quad (x\in \mathcal{X}, B\in \mathcal{B}). 
\end{align*}    
We define the \textbf{conditional mutual information} by
\[  I(X; Y|Z) = \int_{\mathcal{Z}} I(\mu_z, \nu_z) \, d\lambda(z). \]
We have a chain rule \cite[p. 214]{Gray_entropy}
\[  I\left(X;(Y,Z)\right) = I(X;Z) + I(X; Y|Z).  \]
Here the left-hand side is the mutual information between $X$ and $(Y,Z)$.

If $X$ and $Y$ are conditionally independent given $Z$, then $\nu_z(B|x) = \mathbb{P}\left(Y\in B|Z=z\right)$ and hence 
we have $I(\mu_z, \nu_z) = 0$ for $\lambda$-a.e. $z\in \mathcal{Z}$ and $I(X; Y|Z)=0$.

\begin{lemma}[Data-Processing inequality] \label{lemma: data processing inequality}
Let $X, Y, Z$ be random variables that take values in standard Borel spaces.
If $X, Y, Z$ form a Markov chain in this order (i.e. $X$ and $Z$ are conditionally independent given $Y$) then
\[  I(X; Z) \leq I(X; Y). \]
\end{lemma}
  
\begin{proof}
This is a well-known inequality \cite[Theorem 2.8.1]{Cover--Thomas}.
Here we explain a brief proof.
By the chain rule, we have 
\begin{align*}
  I\left(X; (Y, Z)\right) & = I(X; Y) + I(X; Z|Y) \\
   & = I(X;Z) + I(X; Y|Z).
\end{align*}   
Since $X$ and $Z$ are conditionally independent given $Y$, we have $I(X; Z|Y) = 0$.
Hence $I(X; Y) = I(X; Z) + I(X; Z|Y)$.
The conditional mutual information $I(X; Z|Y)$ is nonnegative.
Thus $I(X; Y) \geq I(X;Z)$.
\end{proof}

\begin{proposition}[$I(\mu, \nu)$ is concave in $\mu$ and convex in $\nu$] 
\label{proposition: concavity and convexity of mutual information}
Let $(\mathcal{X}, \mathcal{A})$ and $(\mathcal{Y}, \mathcal{B})$ be standard Borel spaces, and let 
$(\mathcal{Z}, \mathcal{C}, m)$ be a probability space.
  \begin{enumerate}
    \item[(1)] Let $\nu$ be a transition probability on $\mathcal{X}\times \mathcal{Y}$. 
    Suppose that we are given a probability measure $\mu_z$ on $\mathcal{X}$ for each $z\in \mathcal{Z}$ and that 
    $\mu_z$ is measurable in $z$ (namely, the map $\mathcal{Z}\ni z\mapsto \mu_z(A)\in [0,1]$ is measurable for every 
    $A\in \mathcal{A}$). We define a probability measure $\mu$ on $(\mathcal{X}, \mathcal{A})$ by 
    \[  \mu(A) = \int_{\mathcal{Z}} \mu_z(A) \, dm(z), \quad (A\in \mathcal{A}). \]
    Then we have 
    \[  I(\mu, \nu) \geq  \int_{\mathcal{Z}} I(\mu_z, \nu)\, dm(z). \]
    \item[(2)] Let $\mu$ be a probability measure on $\mathcal{X}$. Suppose that we are given a transition probability 
    $\nu_z$ on $\mathcal{X}\times \mathcal{Y}$ for each $z\in \mathcal{Z}$ such that the map 
    $\mathcal{X}\times \mathcal{Z}\ni (x, z) \mapsto \nu_z(B|x)\in [0,1]$ is measurable with respect to $\mathcal{A}\otimes \mathcal{C}$
    for each $B\in \mathcal{B}$.
    We define a transition probability $\nu$ on $\mathcal{X}\times \mathcal{Y}$ by 
    \[  \nu(B|x) = \int_{\mathcal{Z}}\nu_z(B|x)\, dm(z), \quad (x\in \mathcal{X}, B \in \mathcal{B}). \]
    Then we have 
    \[  I(\mu, \nu) \leq  \int_{\mathcal{Z}}I(\mu, \nu_z) \, dm(z). \]
  \end{enumerate}
\end{proposition}

\begin{proof}
See \cite[Proposition 2.10]{Tsukamoto_Variational_R^d}.
In the present paper we will use only the statement (1), but we have also mentioned to (2) for completeness. 
\end{proof}

The following proposition provides a tool for obtaining a lower bound on rate distortion function.

\begin{proposition} \label{proposition: lower bound on mutual information}
Let $\varepsilon >0$ and $a\geq 0$ be real numbers.
Let $\mathcal{X}$ and $\mathcal{Y}$ be measurable spaces with a measurable function 
$\rho\colon \mathcal{X}\times \mathcal{Y}\to [0, \infty)$.
Let $\mu$ be a probability measure on $\mathcal{X}$, and let $\lambda \colon \mathcal{X}\to [0, \infty)$ be a measurable function 
satisfying 
\[  \int_{\mathcal{X}} \lambda(x) 2^{-a\rho(x, y)} d\mu(x) \leq  1 \]
for all $y\in \mathcal{Y}$.
Let $X$ and $Y$ be random variables that take values in $\mathcal{X}$ and $\mathcal{Y}$ respectively and satisfy 
\[  \mathrm{Law} X = \mu, \quad  \mathbb{E} \rho(X, Y) < \varepsilon. \]
Then we have 
\[  I(X; Y) \geq  -a\varepsilon + \int_{\mathcal{X}} \log \lambda(x) \, d\mu(x). \]
\end{proposition}

\begin{proof}
See \cite[Proposition 2.12]{Tsukamoto_Variational_R^d}.
\end{proof}

\subsection{Rate distortion function}  \label{subsection: rate distortion function}

The purpose of this subsection is to introduce rate distortion function for $\mathbb{R}^d$-actions.
Rate distortion function was originally introduced by Shannon \cite{Shannon, Shannon59}.
Readers can find a nice introduction in the book of Cover--Thomas \cite[Chapter 10]{Cover--Thomas}.
Rate distortion theory for continuous-time stochastic processes was investigated by Pursley--Gray \cite{Pursley--Gray}.

We denote by $\mathbf{m}$ the standard Lebesgue measure on $\mathbb{R}^d$.
Let $(\mathcal{X}, \mathbf{d})$ be a compact metric space.
For $L>0$, we define $L^1\left([0,L)^d, \mathcal{X}\right)$ as the space of all measurable maps $f\colon [0,L)^d \to \mathcal{X}$, 
where we identify two maps if they agree $\mathbf{m}$-almost everywhere.
We define a metric $D$ on $L^1\left([0,L)^d, \mathcal{X}\right)$ by 
\[ D(f, g) = \int_{[0,L)^d} \mathbf{d}\left(f(u), g(u)\right) \, d\mathbf{m}(u), \quad \left(f, g\in L^1\left([0,L)^d, \mathcal{X}\right)\right). \]
$L^1\left([0,L)^d, \mathcal{X}\right)$ becomes a complete separable metric space with respect to this metric
\cite[Lemma 2.14]{Tsukamoto_Variational_R^d}. 
Hence it is a standard Borel space with respect to the Borel $\sigma$-algebra.

Let $T\colon \mathbb{R}^d\times \mathcal{X}\to \mathcal{X}$ be a continuous action of $\mathbb{R}^d$ on $\mathcal{X}$.
Let $\mu$ be a $T$-invariant Borel probability measure on $\mathcal{X}$.
Let $\varepsilon>0$ and $L>0$.
We define $R_L(\mathbf{d}, \mu, \varepsilon)$ as the infimum of the mutual information $I(X;Y)$ where $X$ and $Y$ are random variables defined 
on some common probability space $(\Omega, \mathcal{F}, \mathbb{P})$ such that 
  \begin{itemize}
    \item $X$ takes values in $\mathcal{X}$ with $\mathrm{Law} X = \mu$.
    \item $Y$ takes values in $L^1([0,L)^d, \mathcal{X})$ and satisfies 
    \[  \mathbb{E}\left(\frac{1}{L^d} \int_{[0,L)^d} \mathbf{d}\left(T^u X, Y_u\right) \, d\mathbf{m}(u)\right) < \varepsilon, \]  
    where $Y_u = Y_u(\omega)$ $(\omega\in \Omega)$ denotes the value of the function 
    $Y(\omega)\in L^1([0,L)^d, \mathcal{X})$ at $u\in [0,L)^d$.
  \end{itemize}
We define the \textbf{rate distortion function} $R(\mathbf{d}, \mu, \varepsilon)$ by 
\[  R(\mathbf{d}, \mu, \varepsilon) = \lim_{L\to \infty} \frac{R_{L}\left(\mathbf{d}, \mu, \varepsilon \right)}{L^d}. \]
This limit exists and is equal to the infimum of $\frac{R_{L}\left(\mathbf{d}, \mu, \varepsilon \right)}{L^d}$ over $L>0$
\cite[Lemma 2.17]{Tsukamoto_Variational_R^d}.
The value of $R(\mathbf{d}, \mu, \varepsilon)$ is a nonnegative real number.

It is easy to see that \cite[Lemma 2.16 (1)]{Tsukamoto_Variational_R^d} 
\begin{equation} \label{eq: pre-rate distortion function and covering number}
  R_L(\mathbf{d}, \mu, \varepsilon) \leq  \log \#\left(\mathcal{X}, \mathbf{d}_L, \varepsilon\right), 
\end{equation}
where the right-hand side is the logarithm of the $\varepsilon$-covering number 
(introduced in \S \ref{section: main results}) with respect to the metric 
$\mathbf{d}_L(x, y) = \sup_{u\in [0, L)^d} \mathbf{d}(T^u x, T^u y)$.
In particular 
\begin{equation} \label{eq: rate distortion function and varepsilon entropy}
  R(\mathbf{d}, \mu, \varepsilon) \leq  \lim_{L\to \infty} \frac{ \log \#\left(\mathcal{X}, \mathbf{d}_L, \varepsilon\right)}{L^d}
   = S(\mathcal{X}, T, \mathbf{d}, \varepsilon).
\end{equation}
Here $S(\mathcal{X}, T, \mathbf{d}, \varepsilon)$ is the entropy at the scale $\varepsilon$ introduced in 
(\ref{eq: entropy at the scale epsilon}).

We recall that we have defined the upper and lower rate distortion dimensions by 
\[  \urdim(\mathcal{X}, T, \mathbf{d}, \mu) = \limsup_{\varepsilon\to 0} \frac{R(\mathbf{d}, \mu, \varepsilon)}{\log(1/\varepsilon)},\quad 
    \lrdim(\mathcal{X}, T, \mathbf{d}, \mu) = \liminf_{\varepsilon\to 0} \frac{R(\mathbf{d}, \mu, \varepsilon)}{\log(1/\varepsilon)}. \]
These are nonnegative real numbers (possibly $+\infty$). 
When they coincide, we denote the common value by $\rdim\left(\mathcal{X}, T, \mathbf{d}, \mu\right)$.    

The following proposition is an immediate consequence of Proposition \ref{proposition: lower bound on mutual information}.
This will be used in the proof of Theorem \ref{theorem: irregular behavior of upper rate distortion dimension}.

\begin{proposition} \label{proposition: lower bound on rate distortion function}
Let $\varepsilon >0$, $a\geq 0$ and $L>0$ be real numbers.
Suppose that a measurable function $\lambda\colon \mathcal{X}\to [0, \infty)$ satisfies 
\[  \int_{\mathcal{X}} \lambda(x) 2^{-\frac{a}{L^d}\int_{[0, L)^d}\mathbf{d}(T^u x, y_u)d\mathbf{m}(u)} \, d\mu(x) \leq 1 \]
for all $y\in L^1\left([0, L)^d, \mathcal{X}\right)$. Then we have 
\[ R_{L}\left(\mathbf{d}, \mu, \varepsilon\right) \geq - a\varepsilon + \int_{\mathcal{X}} \log \lambda(x) \, d\mu(x). \]
\end{proposition}

\begin{proof}
Apply Proposition \ref{proposition: lower bound on mutual information} to the spaces 
$\mathcal{X}, \mathcal{Y} := L^1\left([0, L)^d, \mathcal{X}\right)$ and the function
\[    \rho(x, y) := \frac{1}{L^d}\int_{[0, L)^d}\mathbf{d}\left(T^u x, y_u\right)\, d\mathbf{m}(u) \quad 
    (x\in \mathcal{X}, y\in L^1\left([0, L)^d, \mathcal{X}\right)). \]
\end{proof}

\begin{remark}
We briefly explain how to modify the definition of rate distortion function for the case of $\mathbb{Z}^d$-actions.
(See also the paper of Huo--Yuan \cite[\S 2.5]{Huo--Yuan}.)
Let $T\colon \mathbb{Z}^d\times \mathcal{X}\to \mathcal{X}$ be a continuous action of the group $\mathbb{Z}^d$
on a compact metric space $(\mathcal{X}, \mathbf{d})$.
Let $\mu$ be a $T$-invariant Borel probability measure on $\mathcal{X}$.
For $\varepsilon>0$ and $L\in \mathbb{N}$ we define $R_L(\mathbf{d}, \mu, \varepsilon)$ as the infimum of 
$I(X; Y)$ where $X$ and $Y$ are random variables defined on 
some common probability space such that 
  \begin{itemize}
    \item $X$ takes values in $\mathcal{X}$ with $\mathrm{Law} X  = \mu$, 
    \item $Y = (Y_{\mathbf{n}})_{\mathbf{n}\in \{0, 1, 2\dots, L-1\}^d}$ and each $Y_{\mathbf{n}}$ takes values in $\mathcal{X}$ with 
    \[  \mathbb{E}\left(\frac{1}{L^d} \sum_{\mathbf{n}\in \{0, 1, 2, \dots, L-1\}^d}
         \mathbf{d}\left(T^{\mathbf{n}}X, Y_{\mathbf{n}}\right)\right) < \varepsilon. \]
  \end{itemize}
We define rate distortion function by 
\[ R(\mathbf{d}, \mu, \varepsilon) = \lim_{L\to \infty} \frac{R_L(\mathbf{d}, \mu, \varepsilon)}{L^d}. \]
This limits always exists and is equal to $\inf_{L\in \mathbb{N}} \frac{R_L(\mathbf{d}, \mu, \varepsilon)}{L^d}$.
Finally we define upper and lower rate distortion dimensions by 
\[  \urdim(\mathcal{X}, T, \mathbf{d}, \mu) = \limsup_{\varepsilon\to 0} \frac{R(\mathbf{d}, \mu, \varepsilon)}{\log(1/\varepsilon)},\quad 
    \lrdim(\mathcal{X}, T, \mathbf{d}, \mu) = \liminf_{\varepsilon\to 0} \frac{R(\mathbf{d}, \mu, \varepsilon)}{\log(1/\varepsilon)}. \]
\end{remark}

\section{Proofs of the main theorems}

The purpose of this section is to prove Theorems \ref{theorem: upper rate distortion dimension and convex combination} and 
\ref{theorem: lower rate distortion dimension and convex combination}.

\subsection{Wasserstein metric and rate distortion function} \label{subsection: Wasserstein metric and rate distortion function}

Here we prepare some basic properties of rate distortion function.
As a technical device for the proofs, we need to use a well-known fact about the Wasserstein metric.

Let $(\mathcal{X}, \mathbf{d})$ be a compact metric space, and let 
$\mathscr{M}(\mathcal{X})$ be the space of Borel probability measures on it.
For $\mu, \nu \in \mathscr{M}(\mathcal{X})$ we define the \textbf{Wasserstein metric} $W_{\mathbf{d}}(\mu, \nu)$ as the infimum of 
\[  \int_{\mathcal{X}\times \mathcal{X}} \mathbf{d}(x, y) \, d\pi(x, y) \]
where $\pi$ runs over Borel probability measures on $\mathcal{X}\times \mathcal{X}$ whose first and second marginals are given by 
$\mu$ and $\nu$ respectively. In other words, $W_{\mathbf{d}}(\mu, \nu)$ is equal to the infimum of 
\[ \mathbb{E} \mathbf{d}(X, Y) \]
where $X$ and $Y$ are random variables defined on a common probability space and taking values in $\mathcal{X}$ with 
$\mathrm{Law} X= \mu$ and $\mathrm{Law} Y = \nu$. Notice that $W_{\mathbf{d}}(\mu, \nu)$ is always finite because 
$(\mathcal{X}, \mathbf{d})$ is compact (and hence, in particular, bounded).

A key fact for us is that the Wasserstein metric metrizes the weak$^*$ topology \cite[Theorem 6.9]{Villani}.
Namely, a sequence $\{\mu_n\}$ in $\mathscr{M}(\mathcal{X})$ converges to $\mu$ in the weak$^*$ topology if and only if 
$W_{\mathbf{d}}(\mu_n, \mu)\to 0$ as $n\to \infty$.

Let $T\colon \mathbb{R}^d\times \mathcal{X}\to \mathcal{X}$ be a continuous action of $\mathbb{R}^d$ on a compact metric space 
$(\mathcal{X}, \mathbf{d})$.
For $L>0$ we define a metric $\bar{\mathbf{d}}_L$ on $\mathcal{X}$ by 
\[ \bar{\mathbf{d}}_L(x, y) = \frac{1}{L^d} \int_{[0, L)^d} \mathbf{d}\left(T^u x, T^u y\right) d\mathbf{m}(u). \]
This metric defines the same topology as the original one of $(\mathcal{X},\mathbf{d})$.
In particular, a sequence $\{\mu_n\}$ in $\mathscr{M}(\mathcal{X})$ converges to $\mu$ in the weak$^*$ topology if and only if 
$W_{\bar{\mathbf{d}}_L}(\mu_n, \mu)\to 0$ as $n\to \infty$.

Recall that, for $\mu\in \mathscr{M}^T(\mathcal{X})$ and $\varepsilon>0$, the function $R_L(\mathbf{d}, \mu, \varepsilon)$
is defined as the infimum of $I(X;Y)$ where $X$ and $Y$ are random variables such that 
$X$ takes values in $\mathcal{X}$ with $\mathrm{Law} X = \mu$ and that 
$Y$ takes values in $L^1\left([0, L)^d, \mathcal{X}\right)$ with 
\[ \mathbb{E}\left(\frac{1}{L^d} \int_{[0, L)^d} \mathbf{d}\left(T^u X, Y_u\right)\, d\mathbf{m}(u)\right) < \varepsilon. \]

\begin{lemma} \label{lemma: convex in varepsilon}
$R_L(\mathbf{d}, \mu, \varepsilon)$ is convex and monotone non-increasing in $\varepsilon$.
In particular, it is a continuous function of the variable $\varepsilon$. 
\end{lemma}

\begin{proof}
The monotonicity is obvious.
We prove that it is convex in $\varepsilon$.
Let $\varepsilon_0$ and $\varepsilon_1$ be positive numbers.
Let $p_0, p_1\in [0,1]$ with $p_0+p_1=1$.
For any $\delta>0$ there exist random variables $(X^{(i)}, Y^{(i)})$ $(i=0,1)$ that take values in 
$\mathcal{X}\times L^1([0, L)^d, \mathcal{X})$ with the following conditions:
 \begin{itemize}
    \item 
    $\mathrm{Law} X^{(i)} = \mu$ and 
    $\mathbb{E}\left(\frac{1}{L^d}\int_{[0, L)^d}\mathbf{d}\left(T^t X^{(i)}, Y^{(i)}_t\right) d\mathbf{m}(t)\right)< \varepsilon_i$,
    \item $I\left(X^{(i)}; Y^{(i)}\right) < R_L(\mathbf{d}, \mu, \varepsilon_i) + \delta$, 
    \item $(X^{(0)}, Y^{(0)})$ and $(X^{(1)}, Y^{(1)})$ are independent.
 \end{itemize}
Moreover we can take a random variable $Z$ independent of $(X^{(0)}, Y^{(0)}, X^{(1)}, Y^{(1)})$ such that 
\[  Z = \begin{cases} 0  & \text{in probability $p_0$} \\ 1 & \text{in probability $p_1$} \end{cases}. \]
We set $X=X^{(Z)}$ and $Y = Y^{(Z)}$.
Since $X^{(0)}$ and $X^{(1)}$ have the same distribution, $X$ and $Z$ are independent.
Hence $I(X; Z) = 0$.
We have 
\begin{align*}
  \mathbb{E}\left(\frac{1}{L^d}\int_{[0, L)^d} \mathbf{d}\left(T^t X, Y_t\right) \, d\mathbf{m}(t)\right) 
  & = \sum_{i=0, 1} p_i \mathbb{E}\left(\frac{1}{L^d}\int_{[0, L)^d} \mathbf{d}\left(T^t X^{(i)}, Y^{(i)}_t\right) \, d\mathbf{m}(t)\right) \\
  & < p_0 \varepsilon_0 + p_1 \varepsilon_1.
\end{align*}
We also have 
\begin{align*}
  I(X; Y) & \leq I\left(X; (Y, Z)\right) = \underbrace{I(X; Z)}_{=0} + I(X;Y|Z) \\
  & = p_0 I\left(X^{(0)}; Y^{(0)}\right) + p_1 I\left(X^{(1)}; Y^{(1)}\right)  \\
  & < p_0 R_L(\mathbf{d}, \mu, \varepsilon_0) + p_1 R_L(\mathbf{d}, \mu, \varepsilon_1) + \delta.
\end{align*}
This implies 
\[  R_L(\mathbf{d}, \mu, p_0 \varepsilon_0 + p_1 \varepsilon_1) <
 p_0 R_L(\mathbf{d}, \mu, \varepsilon_0) + p_1 R_L(\mathbf{d}, \mu, \varepsilon_1) + \delta. \]
Letting $\delta\to 0$, we get 
\[ R_L(\mathbf{d}, \mu, p_0 \varepsilon_0 + p_1 \varepsilon_1) \leq
 p_0 R_L(\mathbf{d}, \mu, \varepsilon_0) + p_1 R_L(\mathbf{d}, \mu, \varepsilon_1). \]
\end{proof}

\begin{lemma} \label{lemma: continuity of pre-rate distortion function}
$R_L(\mathbf{d}, \mu, \varepsilon)$ is continuous in $(\mu, \varepsilon)$.
Namely, if $\varepsilon_n \to \varepsilon$ and 
if $\mu_n\to \mu$ in $\mathscr{M}^T(\mathcal{X})$ with respect to the weak$^*$ topology,
then
$R_L(\mathbf{d}, \mu_n, \varepsilon_n)\to R_L(\mathbf{d}, \mu, \varepsilon)$.
\end{lemma}

\begin{proof}
We first show $\limsup_{n\to \infty} R_L(\mathbf{d}, \mu_n, \varepsilon_n)\leq  R_L(\mathbf{d}, \mu,\varepsilon)$.
Let $\delta>0$ be arbitrary.
There exist random variables $X$ and $Y$ that take values in $\mathcal{X}$ and $L^1([0, L)^d, \mathcal{X})$ 
respectively and satisfy
\[ \mathrm{Law} X = \mu, \quad 
    \mathbb{E}\left(\frac{1}{L^d} \int_{[0, L)^d} \mathbf{d}\left(T^t X, Y_t\right)\, d\mathbf{m}(t)\right) < \varepsilon, \quad 
    I(X; Y) < R_L(\mathbf{d}, \mu, \varepsilon) + \delta. \]
Take $\eta>0$ satisfying 
\[ \mathbb{E}\left(\frac{1}{L^d} \int_{[0, L)^d} \mathbf{d}\left(T^t X, Y_t\right)\, d\mathbf{m}(t)\right) + \eta < \varepsilon_n \]
for all sufficiently large $n$.
We have $W_{\bar{\mathbf{d}}_L}(\mu_n, \mu) < \eta$ for large $n$.
Hence we can assume that (for any large $n$) there is a random variable $X_n$ that takes values in $\mathcal{X}$ 
and is conditionally independent of $Y$ given $X$ and satisfies
\[ \mathrm{Law} X_n = \mu_n, \quad  \mathbb{E}\bar{\mathbf{d}}_L(X_n, X) < \eta. \]
Then 
\begin{align*}
     \mathbb{E}\left(\frac{1}{L^d} \int_{[0, L)^d} \mathbf{d}\left(T^t X_n, Y_t\right)\, d\mathbf{m}(t)\right)
     &  \leq   \mathbb{E}\bar{\mathbf{d}}_L(X_n, X) 
     + \mathbb{E}\left(\frac{1}{L^d} \int_{[0, L)^d} \mathbf{d}\left(T^t X, Y_t\right)\, d\mathbf{m}(t)\right) \\
     & < \eta + \mathbb{E}\left(\frac{1}{L^d} \int_{[0, L)^d} \mathbf{d}\left(T^t X, Y_t\right)\, d\mathbf{m}(t)\right) \\
      & < \varepsilon_n.
\end{align*}      
The random variables $Y, X, X_n$ form a Markov chain in this order. Hence by the data-processing inequality 
\[  I(X_n;Y) \leq I(X; Y) <  R_L(\mathbf{d}, \mu, \varepsilon) + \delta. \]
Hence $R_L(\mathbf{d}, \mu_n, \varepsilon_n) < R_L(\mathbf{d}, \mu, \varepsilon)+\delta$ for large $n$.
In particular, $\limsup_{n\to \infty} R_L(\mathbf{d}, \mu_n, \varepsilon_n)\leq  R_L(\mathbf{d}, \mu,\varepsilon) + \delta$.
Since $\delta>0$ is arbitrary, this shows $\limsup_{n\to \infty} R_L(\mathbf{d}, \mu_n, \varepsilon_n)\leq  R_L(\mathbf{d}, \mu,\varepsilon)$.

Next we prove $\liminf_{n\to \infty} R_L(\mathbf{d}, \mu_n, \varepsilon_n)\geq  R_L(\mathbf{d}, \mu,\varepsilon)$.
The argument is similar to the above.
Let $\delta>0$ be arbitrary.
We take a sufficiently large $n\geq 1$ such that $\varepsilon_n < \varepsilon + \delta$ and 
$W_{\bar{\mathbf{d}}_L}(\mu_n, \mu) < \delta$.
For such $n$, we take random variables $X_n$ and $Y$ that take values in $\mathcal{X}$ and 
$L^1([0, L)^d, \mathcal{X})$ respectively and satisfy
\[ \mathrm{Law} X_n = \mu_n, \quad 
    \mathbb{E}\left(\frac{1}{L^d}\int_{[0, L)^d} \mathbf{d}\left(T^t X_n, Y_t\right) \, d\mathbf{m}(t)\right) < \varepsilon_n, \quad 
    I(X_n;Y) < R_L(\mathbf{d}, \mu_n, \varepsilon_n) + \delta. \]
Since $W_{\bar{\mathbf{d}}_L}(\mu_n, \mu) < \delta$, 
we can also assume that there is a random variable $X$ that takes values in $\mathcal{X}$
with $\mathrm{Law} X = \mu$ and $\mathbb{E} \bar{\mathbf{d}}_L (X,X_n) < \delta$ and 
is conditionally independent of $Y$ given $X_n$. 
Then we have 
\[   \mathbb{E}\left(\frac{1}{L^d}\int_{[0, L)^d} \mathbf{d}\left(T^t X, Y_t\right) \, d\mathbf{m}(t)\right) < \varepsilon_n + \delta 
      < \varepsilon + 2\delta. \]
Here we have used $\varepsilon_n < \varepsilon + \delta$.
By the data-processing inequality 
\[  I(X; Y) \leq  I(X_n; Y) < R_L(\mathbf{d}, \mu_n, \varepsilon_n) + \delta. \]
Hence $R_L(\mathbf{d}, \mu, \varepsilon+2\delta) \leq R_L(\mathbf{d}, \mu_n, \varepsilon_n)+\delta$.
Therefore we have $R_L(\mathbf{d}, \mu, \varepsilon+2\delta) \leq \liminf_{n\to \infty} R_L(\mathbf{d}, \mu_n, \varepsilon_n) + \delta$.
By Lemma \ref{lemma: convex in varepsilon}
\[  R_L(\mathbf{d}, \mu, \varepsilon) = \lim_{\delta\to 0} R_L(\mathbf{d}, \mu, \varepsilon+2\delta). \]
Thus we conclude that 
$R_L(\mathbf{d}, \mu, \varepsilon) \leq  \liminf_{n\to \infty} R_L(\mathbf{d}, \mu_n, \varepsilon_n)$.
\end{proof}

Recall that the rate distortion function $R(\mathbf{d}, \mu, \varepsilon)$ $(\mu\in \mathscr{M}^T(\mathcal{X}), \varepsilon>0)$ is defined by 
\[ R(\mathbf{d}, \mu, \varepsilon) = \lim_{L\to \infty} \frac{R_L(\mathbf{d}, \mu, \varepsilon)}{L^d}
    = \inf_{L>0} \frac{R_L(\mathbf{d},\mu,\varepsilon)}{L^d}. \]

\begin{lemma}   \label{lemma: upper semi-continuity of rate distortion function}
The rate distortion function $R(\mathbf{d}, \mu, \varepsilon)$ is upper semi-continuous in $(\mu, \varepsilon)$.
Namely, if $\varepsilon_n \to \varepsilon$ and 
if $\mu_n\to \mu$ in $\mathscr{M}^T(\mathcal{X})$ with respect to the weak$^*$ topology,
then
\[  \limsup_{n\to \infty} R(\mathbf{d}, \mu_n, \varepsilon_n) \leq  R(\mathbf{d}, \mu, \varepsilon). \]
In particular, $R(\mathbf{d}, \mu, \varepsilon)$ is a measurable function of $(\mu, \varepsilon)$.
\end{lemma}

\begin{proof}
Let $\delta>0$ be arbitrary.
There is $L>0$ such that $R(\mathbf{d}, \mu, \varepsilon) +\delta > \frac{1}{L^d}R_L(\mathbf{d}, \mu, \varepsilon)$.
By Lemma \ref{lemma: continuity of pre-rate distortion function}, for all sufficiently large $n$, we have 
$R_L(\mathbf{d}, \mu_n, \varepsilon_n) < R_L(\mathbf{d}, \mu, \varepsilon) + \delta$.
Hence 
\[ R(\mathbf{d},\mu_n,\varepsilon_n) \leq 
   \frac{1}{L^d}R_L(\mathbf{d}, \mu_n,\varepsilon_n) 
   < R(\mathbf{d}, \mu, \varepsilon)+\delta+\frac{\delta}{L^d}. \]
Therefore $\limsup_{n\to \infty} R(\mathbf{d},\mu_n,\varepsilon_n) \leq R(\mathbf{d}, \mu, \varepsilon) + \delta+\frac{\delta}{L^d}$.
Since $\delta>0$ is arbitrary, this shows the lemma.
\end{proof}

\begin{remark}
We will see below (Lemma \ref{lemma: convexity in finite sum}) that $R(\mathbf{d}, \mu, \varepsilon)$ is continuous in $\varepsilon$.
However it may not be continuous in $\mu$ in general.
For the case of $\mathbb{Z}$-actions, a simple example is given as follows.
Let $\mathcal{X} = [0,1]^{\mathbb{Z}}$ be the infinite dimensional cube with a metric 
$\mathbf{d}(x, y) := \sup_{n\in \mathbb{Z}}2^{-|n|}|x_n-y_n|$ and 
the natural shift map $T\colon \mathcal{X}\to \mathcal{X}$.
For $n\geq 1$, let $K_n = \{x\in \mathcal{X}\mid T^n x = x \}$.
There is a natural homeomorphism $\varphi_n\colon [0,1]^{n}\to K_n$ defined by 
\[ \varphi_n(x_0, x_1, \dots, x_{n-1})_{\ell+m n} =x_{\ell} \quad (0\leq \ell \leq n-1,  m\in \mathbb{Z}). \] 
We denote by $\mathbf{m}_n$ the Lebesgue measure on $[0,1]^n$ and set $\nu_n = (\varphi_n)_*\mathbf{m}_n$.
Define $\mu_n = \frac{1}{n}\sum_{k=0}^{n-1} T^k_*\nu_n \in \mathscr{M}^T(\mathcal{X})$.
Then the sequence $\mu_n$ converges to $\mu:=\mathbf{m}_1^{\otimes \mathbb{Z}}$.
On the one hand we have 
\[  R(\mathbf{d}, \mu, \varepsilon) \sim \log (1/\varepsilon). \]
On the other hand, for every $n\in \mathbb{Z}$ we have 
\[  R(\mathbf{d}, \mu_n, \varepsilon) = 0. \]
Therefore the rate distortion function is not continuous in $\mu$.
A similar example can be constructed for $\mathbb{R}^d$-actions by using the method of \S \ref{section: construction of an example}.
\end{remark}

\subsection{Rate distortion function and convex combination of measures} 
\label{subsection: rate distortion function and convex combination of measures}

Here we prove the formula (\ref{eq: rate distortion function and convex combination}) in \S \ref{section: main results}.
As we mentioned, this formula is well-known for $\mathbb{Z}$-actions in the classical information theory literature 
\cite{Gray--Davisson, Shields--Neuhoff--Davisson--Ledrappier,
LeonGarcia--Davisson--Neuhoff, Effros--Chou--Gray}.
The proof consists of several intermediate lemmas.
All the arguments (including the ones in \S \ref{subsection: Wasserstein metric and rate distortion function}) 
are basic and standard.
Hence it is hopefully easy for motivated readers to modify the arguments 
for $\mathbb{Z}^d$-actions if they have interests in the case of $\mathbb{Z}^d$.

Throughout of this subsection we assume that $(\mathcal{X}, \mathbf{d})$ is a compact metric space with 
a continuous action $T\colon \mathbb{R}^d\times \mathcal{X}\to \mathcal{X}$.

\begin{lemma} \label{lemma: convexity in finite sum}
Let $p_0, p_1 \in [0,1]$ be real numbers with $p_0+p_1=1$.
Let $\mu_0, \mu_1\in \mathscr{M}^T(\mathcal{X})$.
For positive numbers $\varepsilon_0$ and $\varepsilon_1$, we have 
\[   R(\mathbf{d}, p_0\mu_0+p_1\mu_1, p_0\varepsilon_0+p_1\varepsilon_1) \leq 
      p_0 R(\mathbf{d}, \mu_0, \varepsilon_0)+p_1 R(\mathbf{d}, \mu_1, \varepsilon_1). \]
In particular, the rate distortion function $R(\mathbf{d}, \mu, \varepsilon)$ is a continuous function of the variable $\varepsilon$.
\end{lemma}

\begin{proof}
The proof is very close to Lemma \ref{lemma: convex in varepsilon}.
Let $L>0$ and $\delta>0$.
There exist random variables $(X^{(i)}, Y^{(i)})$ $(i=0,1)$ that take values in 
$\mathcal{X}\times L^1([0, L)^d, \mathcal{X})$ with the following properties:
 \begin{itemize}
    \item 
    $\mathrm{Law} X^{(i)} = \mu_i$ and 
    $\mathbb{E}\left(\frac{1}{L^d}\int_{[0, L)^d}\mathbf{d}\left(T^t X^{(i)}, Y^{(i)}_t\right) d\mathbf{m}(t)\right)< \varepsilon_i$,
    \item $I\left(X^{(i)}; Y^{(i)}\right) < R_L(\mathbf{d}, \mu, \varepsilon_i) + \delta$, 
    \item $(X^{(0)}, Y^{(0)})$ and $(X^{(1)}, Y^{(1)})$ are independent.
 \end{itemize}
Moreover we can take a random variable $Z$ independent of $(X^{(0)}, Y^{(0)}, X^{(1)}, Y^{(1)})$ such that 
\[  Z = \begin{cases} 0  & \text{in probability $p_0$} \\ 1 & \text{in probability $p_1$} \end{cases}. \]
We set $X=X^{(Z)}$ and $Y = Y^{(Z)}$. 
(So far the argument is the same with the proof of Lemma 
\ref{lemma: convex in varepsilon}. The difference is that $X$ and $Z$ are not necessarily independent here.
However this is very minor.)

We have $\mathrm{Law} X = p_0\mu_0+p_1\mu_1$ and 
\[ \mathbb{E}\left(\frac{1}{L^d}\int_{[0, L)^d}\mathbf{d}\left(T^t X, Y_t\right) \, d\mathbf{m}(t)\right) < p_0 \varepsilon_0 + p_1 \varepsilon_1. \]
As in the proof of Lemma \ref{lemma: convex in varepsilon} we also have
\begin{align*}
  I(X; Y) & \leq I\left(X; (Y, Z)\right) \\   
  & = \underbrace{I(X; Z)}_{\leq H(Z) \leq 1} + I(X; Y|Z) \\
  & \leq 1 + p_0 R_L(\mathbf{d}, \mu_0, \varepsilon_0) + p_1 R_L(\mathbf{d}, \mu_1, \varepsilon_1) +\delta.
\end{align*}
Therefore 
\[  R_L(\mathbf{d}, p_0\mu_0+p_1\mu_1, p_0\varepsilon_0+p_1\varepsilon_1) \leq 
    1 + p_0 R_L(\mathbf{d}, \mu_0, \varepsilon_0) + p_1 R_L(\mathbf{d}, \mu_1, \varepsilon_1) +\delta. \]
Dividing this by $L^d$ and letting $L\to \infty$, we get the conclusion.
\end{proof}

Let $\lambda$ be a Borel probability measure on $\mathscr{M}^T(\mathcal{X})$ and define 
\[  \mu = \int_{\mathscr{M}^T(\mathcal{X})} \nu \, d\lambda(\nu) \in \mathscr{M}^T(\mathcal{X}). \]

\begin{lemma} \label{lemma: lower bound of rate distortion function and convex combination of measures}
For a positive number $\varepsilon$
\[    R(\mathbf{d}, \mu, \varepsilon)
   \geq  \inf\left\{\int_{\mathscr{M}^T(\mathcal{X})} R(\mathbf{d},\nu, \varepsilon_\nu)d\lambda(\nu) \middle|\, 
\text{\parbox{7cm}{$\mathscr{M}^T(\mathcal{X})\ni \nu \mapsto \varepsilon_\nu \in (0, \infty)$ with  
$\int_{\mathscr{M}^T(\mathcal{X})} \varepsilon_\nu \, d\lambda(\nu) \leq \varepsilon$}}\right\}. \]
Here the infimum is taken over all measurable functions $\mathscr{M}^T(\mathcal{X}) \ni \nu \mapsto \varepsilon_\nu \in (0, \infty)$
satisfying $\int_{\mathscr{M}^T(\mathcal{X})} \varepsilon_{\nu} \, d\lambda(\nu) \leq \varepsilon$.
\end{lemma}

\begin{proof}
Let $L>0$.
Suppose that random variables $X$ and $Y$ take values in $\mathcal{X}$ and $L^1([0, L)^d, \mathcal{X})$
respectively and satisfy $\mathrm{Law} X = \mu$ and 
$\mathbb{E}\left(\frac{1}{L^d}\int_{[0, L)^d} \mathbf{d}\left(T^t X, Y_t\right)\, d\mathbf{m}(t)\right) < \varepsilon$.
Take $\delta>0$ satisfying 
$\mathbb{E}\left(\frac{1}{L^d}\int_{[0, L)^d} \mathbf{d}\left(T^t X, Y_t\right)\, d\mathbf{m}(t)\right) + \delta < \varepsilon$.

Let $p$ be the regular conditional distribution of $Y$ given $X$.
Namely 
\[  p(A|x) = \mathbb{P}(Y\in A|X=x), \quad (x\in \mathcal{X}, A\subset L^1([0, L)^d, \mathcal{X})).  \]
For each $\nu\in \mathscr{M}^T(\mathcal{X})$ we take random variables $X^{(\nu)}$ and $Y^{(\nu)}$ that take values in 
$\mathcal{X}$ and $L^1([0, L)^d, \mathcal{X})$ respectively and satisfy
$\mathrm{Law} X^{(\nu)} = \nu$ and 
\[ \mathbb{P}(Y^{(\nu)}\in A | X^{(\nu)} = x) = p(A|x), \quad (x\in \mathcal{X}, A\subset L^1([0, L)^d, \mathcal{X})).  \]

We have $I(X; Y) = I(\mu, p) = I\left(\int \nu d\lambda(\nu), p\right)$.
By the concavity of mutual information (Proposition \ref{proposition: concavity and convexity of mutual information} (1))
\[ I(X; Y) \geq \int_{\mathscr{M}^T(\mathcal{X})} I(\nu, p)\, d\lambda(\nu) 
     = \int_{\mathscr{M}^T(\mathcal{X})} I\left(X^{(\nu)}; Y^{(\nu)}\right)  d\lambda(\nu). \]
On the other hand
\begin{align*}
  & \mathbb{E}\left(\int_{[0, L)^d} \mathbf{d}\left(T^t X, Y_t\right) d\mathbf{m}(t)\right) \\
  & = 
  \int_{\mathcal{X}\times L^1([0, L)^d, \mathcal{X})}\left(\int_{[0, L)^d} \mathbf{d}(T^t x, y_t)\, d\mathbf{m}(t) \right)  d\left(\mu(x)p(y|x)\right) \\
  & = \int_{\mathscr{M}^T(\mathcal{X})}
  \left\{\int_{\mathcal{X}\times L^1([0, L)^d, \mathcal{X})}\left(\int_{[0, L)^d} \mathbf{d}(T^t x, y_t)\, d\mathbf{m}(t)\right) 
   d\left(\nu(x)p(y|x)\right)\right\}
   d\lambda(\nu) \\
  & = \int_{\mathscr{M}^T(\mathcal{X})} \mathbb{E}\left(\int_{[0, L)^d} \mathbf{d}\left(T^t X^{(\nu)}, Y^{(\nu)}_t\right) d\mathbf{m}(t)\right)
        d\lambda(\nu). 
\end{align*}
We set 
\[  \varepsilon_{\nu} = \mathbb{E}\left(\frac{1}{L^d} \int_{[0, L)^d} \mathbf{d}\left(T^t X^{(\nu)}, Y^{(\nu)}_t\right)d\mathbf{m}(t)\right) +\delta. \]
We have 
\[ \int_{\mathscr{M}^T(\mathcal{X})} \varepsilon_{\nu}\, d\lambda(\nu) 
  = \mathbb{E}\left(\frac{1}{L^d} \int_{[0, L)^d} \mathbf{d}\left(T^t X, Y_t\right) d\mathbf{m}(t)\right) +\delta < \varepsilon. \]
From the definition of rate distortion function,
\[  \frac{1}{L^d} I\left(X^{(\nu)}; Y^{(\nu)}\right) \geq \frac{1}{L^d} R_L\left(\mathbf{d}, \nu, \varepsilon_{\nu}\right) \geq 
      R\left(\mathbf{d}, \nu, \varepsilon_{\nu}\right). \]
Therefore 
\[  \frac{1}{L^d} I(X; Y) \geq   \frac{1}{L^d} \int_{\mathscr{M}^T(\mathcal{X})} I\left(X^{(\nu)}; Y^{(\nu)}\right)  d\lambda(\nu)
    \geq \int_{\mathscr{M}^T(\mathcal{X})} R\left(\mathbf{d}, \nu, \varepsilon_{\nu}\right)d\lambda(\nu). \]
This proves the claim.
\end{proof}

\begin{lemma} \label{lemma: towards upper bound of rate distortion function and convex combination of measures}
Let $\varepsilon>0$ and let $\mathscr{M}^T(\mathcal{X}) \ni \nu \mapsto \varepsilon_\nu \in (0, \infty)$ be a measurable function
with $\int_{\mathscr{M}^T(\mathcal{X})} \varepsilon_{\nu} \, d\lambda(\nu) \leq \varepsilon$.
We suppose that there is a positive number $c$ satisfying $\varepsilon_{\nu} \geq  c$ for all 
$\nu \in \mathscr{M}^T(\mathcal{X})$. Then 
\[  R(\mathbf{d}, \mu, \varepsilon) \leq \int_{\mathscr{M}^T(\mathcal{X})} R(\mathbf{d}, \nu, \varepsilon_{\nu})\, d\lambda(\nu). \]
\end{lemma}

\begin{proof}
We assume that the space $\mathscr{M}^T(\mathcal{X})$ is endowed with the Wasserstein metric $W_{\mathbf{d}}$.
(Indeed any metric will do the same work.)
For a finite partition $\alpha = \{A_1, \dots, A_n\}$ of $\mathscr{M}^T(\mathcal{X})$, we denote by $\mathrm{mesh}(\alpha)$
the maximum of the diameter $\mathrm{diam} (A_i, W_{\mathbf{d}})$ $(1\leq i \leq n)$.
For two finite partitions $\alpha = \{A_1, \dots, A_n\}$ and $\beta  = \{B_1, \dots, B_m\}$ of 
$\mathscr{M}^T(\mathcal{X})$, we say that $\beta$ is a refinement of $\alpha$ (denoted by $\alpha \prec \beta$)
if every $B_j$ $(1\leq j \leq m)$ is contained in some $A_i$ $(1\leq i \leq n)$.

From (\ref{eq: rate distortion function and varepsilon entropy}) in \S \ref{subsection: rate distortion function}
we have 
$R(\mathbf{d}, \theta, \varepsilon_{\nu}) \leq  S(\mathcal{X}, T, \mathbf{d}, c) < \infty$ for all 
$\nu, \theta \in \mathscr{M}^T(\mathcal{X})$.
We divide the proof into two steps.

\textbf{Step 1:}
We suppose that the function $\mathscr{M}^T(\mathcal{X}) \ni \nu \mapsto \varepsilon_{\nu}$ 
is a simple function. 
Namely, we suppose that there is a finite measurable partition 
$\alpha = \{A_1, \dots, A_n\}$ of $\mathscr{M}^T(\mathcal{X})$ such that $\varepsilon_{\nu}$ is constant on each $A_i$.

We take a sequence of finite measurable partitions $\alpha \prec \beta_1 \prec \beta_2 \prec \dots$ of 
$\mathscr{M}^T(\mathcal{X})$ such that $\mathrm{mesh}(\beta_k)\to 0$ as $k \to \infty$. 
We define a map $\mathcal{\beta}_k\colon \mathscr{M}^T(\mathcal{X})\to \mathscr{M}^T(\mathcal{X})$ as follows.
Let $B\in \mathcal{\beta}_k$.
 \begin{itemize}
    \item If $\lambda(B)>0$ then for any $\nu\in B$ we set 
    \[  \beta_k(\nu) = \frac{1}{\lambda(B)} \int_{B} \theta \, d\lambda(\theta). \]
    \item If $\lambda(B)=0$ then we pick $\theta_B \in B$ and set $\beta_k(\nu) = \theta_B$ for all $\nu \in B$.
 \end{itemize}
We have 
\[ \mu = \int_{\mathscr{M}^T(\mathcal{X})} \nu \, d\lambda(\nu)  
   = \sum_{\substack{B\in \beta_k \\ \lambda(B)>0}}  \lambda(B) \left(\frac{1}{\lambda(B)} \int_B \nu\, d\lambda(\nu)\right)
   = \int_{\mathscr{M}^T(\mathcal{X})} \beta_k(\nu)\, d\lambda(\nu). \]
The function $\varepsilon_{\nu}$ is constant on each $B\in \beta_k$. 
Therefore by Lemma \ref{lemma: convexity in finite sum}
\[ R(\mathbf{d}, \mu, \varepsilon) \leq R\left(\mathbf{d}, \mu, \int \varepsilon_{\nu} d\lambda(\nu)\right)
   \leq \int_{\mathscr{M}^T(\mathcal{X})} R\left(\mathbf{d}, \beta_k(\nu), \varepsilon_{\nu}\right) d\lambda(\nu). \]

For each $\nu \in \mathscr{M}^T(\mathcal{X})$, the sequence $\{\beta_k(\nu)\}_{k=1}^\infty$ converges to
$\nu$ in the weak$^*$ topology since we have assumed $\mathrm{mesh}(\beta_k) \to 0$.
By the upper semi-continuity of rate distortion function (Lemma \ref{lemma: upper semi-continuity of rate distortion function}),
we have $\limsup_{k\to \infty} R\left(\mathbf{d},\beta_k(\nu),\varepsilon_{\nu}\right)\leq R(\mathbf{d},\nu,\varepsilon_{\nu})$.
Recall that we have the uniform bound 
$R\left(\mathbf{d}, \beta_k(\nu), \varepsilon_{\nu}\right) \leq  S(\mathcal{X}, T, \mathbf{d}, c)<\infty$.
Therefore we can use Fatou’s lemma and get 
\begin{align*}
  \limsup_{k\to \infty} \int_{\mathscr{M}^T(\mathcal{X})} R\left(\mathbf{d}, \beta_k(\nu), \varepsilon_{\nu}\right) d\lambda(\nu)
  & \leq \int_{\mathscr{M}^T(\mathcal{X})} \limsup_{k\to \infty} R\left(\mathbf{d}, \beta_k(\nu), \varepsilon_{\nu}\right) d\lambda(\nu) \\
  & \leq  \int_{\mathscr{M}^T(\mathcal{X})} R(\mathbf{d}, \nu, \varepsilon_{\nu})\, d\lambda(\nu).
\end{align*}
Thus $R(\mathbf{d}, \mu, \varepsilon) \leq  \int_{\mathscr{M}^T(\mathcal{X})} R(\mathbf{d}, \nu, \varepsilon_{\nu})\, d\lambda(\nu)$.

\textbf{Step 2:}
We consider the general case. (Namely $\varepsilon_{\nu}$ is not necessarily a simple function.)
We can choose a sequence of simple functions 
$c\leq \varepsilon^{(1)}_\nu \leq \varepsilon^{(2)}_\nu  \leq \varepsilon^{(3)}_{\nu} \leq \dots$
such that we have pointwise convergence $\varepsilon^{(n)}_\nu \to \varepsilon_{\nu}$ as $n \to \infty$ 
at each $\nu \in \mathscr{M}^T(\mathcal{X})$.

We have 
\[  \int_{\mathscr{M}^T(\mathcal{X})} \varepsilon^{(n)}_\nu d\lambda(\nu) \leq 
\int_{\mathscr{M}^T(\mathcal{X})} \varepsilon_\nu d\lambda(\nu) \leq \varepsilon. \] 
We apply the result of Step 1 to $\varepsilon^{(n)}_\nu$ and get 
\[  R(\mathbf{d}, \mu, \varepsilon) \leq  
     \int_{\mathscr{M}^T(\mathcal{X})} R\left(\mathbf{d}, \nu, \varepsilon^{(n)}_\nu\right) d\lambda(\nu). \]
By Lemma \ref{lemma: convexity in finite sum}, $R\left(\mathbf{d}, \nu, \varepsilon^{(n)}_\nu\right) \to R(\mathbf{d}, \nu, \varepsilon_{\nu})$
as $n\to \infty$ at each $\nu \in \mathscr{M}^T(\mathcal{X})$.
Since we have the uniform bound $R\left(\mathbf{d}, \nu, \varepsilon^{(n)}_\nu\right) \leq S(\mathcal{X}, T, \mathbf{d}, c) < \infty$,
we can use Lebesgue’s dominated convergence theorem and get 
\[  \lim_{n\to \infty} \int_{\mathscr{M}^T(\mathcal{X})} R\left(\mathbf{d}, \nu, \varepsilon^{(n)}_\nu\right) d\lambda(\nu)
    =  \int_{\mathscr{M}^T(\mathcal{X})} R(\mathbf{d}, \nu, \varepsilon_{\nu}) \, d\lambda(\nu). \]
Therefore we conclude 
\[ R(\mathbf{d}, \mu, \varepsilon)  \leq  \int_{\mathscr{M}^T(\mathcal{X})} R(\mathbf{d}, \nu, \varepsilon_{\nu}) \, d\lambda(\nu). \]
\end{proof}

\begin{lemma} \label{lemma: upper bound of rate distortion function and convex combination of measures}
Let $\varepsilon>0$ and let $\mathscr{M}^T(\mathcal{X}) \ni \nu \mapsto \varepsilon_\nu \in (0, \infty)$ be a measurable function
with $\int_{\mathscr{M}^T(\mathcal{X})} \varepsilon_{\nu} \, d\lambda(\nu) \leq \varepsilon$.
Then we have 
\[ R(\mathbf{d}, \mu, \varepsilon) \leq \int_{\mathscr{M}^T(\mathcal{X})} R(\mathbf{d}, \nu, \varepsilon_{\nu})\, d\lambda(\nu). \]
\end{lemma}

\begin{proof}
We can assume that $R(\mathbf{d}, \nu, \varepsilon_{\nu})$ is a $\lambda$-integralable function of $\nu$.
(Otherwise the statement is obvious.)
Let $n\geq 1$.
We apply Lemma \ref{lemma: towards upper bound of rate distortion function and convex combination of measures}
to the measurable function $\mathscr{M}^T(\mathcal{X})\ni \nu \mapsto \varepsilon_{\nu} +\frac{1}{n} \in (0, \infty)$.
Notice that this function is bounded from below by $\frac{1}{n}$.
Then we have 
\[ R\left(\mathbf{d}, \mu, \varepsilon+\frac{1}{n}\right) \leq  \int_{\mathscr{M}^T(\mathcal{X})} 
    R\left(\mathbf{d}, \nu, \varepsilon_{\nu}+\frac{1}{n}\right) d\lambda(\nu). \]
By Lemma \ref{lemma: convexity in finite sum}
\[ R(\mathbf{d}, \mu, \varepsilon) = \lim_{n\to \infty} R\left(\mathbf{d}, \mu, \varepsilon+\frac{1}{n}\right). \]
For each $\nu\in \mathscr{M}^T(\mathcal{X})$ we also have 
\[ R\left(\mathbf{d}, \nu, \varepsilon_{\nu}\right) = \lim_{n\to \infty} R\left(\mathbf{d}, \nu, \varepsilon_{\nu}+\frac{1}{n}\right). \]
We notice that $R\left(\mathbf{d}, \nu, \varepsilon_{\nu}+\frac{1}{n}\right) \leq R\left(\mathbf{d}, \nu, \varepsilon_{\nu}\right)$ and that 
the right-hand side is assumed to be $\lambda$-integrable.
Therefore we can use Lebesgue’s dominated convergence theorem and get 
\[ \lim_{n\to \infty} \int_{\mathscr{M}^T(\mathcal{X})} R\left(\mathbf{d}, \nu, \varepsilon_{\nu}+\frac{1}{n}\right) d\lambda(\nu)
    = \int_{\mathscr{M}^T(\mathcal{X})} R\left(\mathbf{d}, \nu, \varepsilon_{\nu}\right) d\lambda(\nu). \]
Thus we conclude that 
\[  R(\mathbf{d}, \mu, \varepsilon) \leq \int_{\mathscr{M}^T(\mathcal{X})} R\left(\mathbf{d}, \nu, \varepsilon_{\nu}\right) d\lambda(\nu). \]
\end{proof}

By combining Lemmas \ref{lemma: lower bound of rate distortion function and convex combination of measures} and
\ref{lemma: upper bound of rate distortion function and convex combination of measures}, we get the conclusion of this 
subsection.

\begin{theorem}  \label{theorem: rate distortion function and convex combination of measures}
Let $\lambda$ be a Borel probability measure on $\mathscr{M}^T(\mathcal{X})$ and define a measure 
$\mu \in \mathscr{M}^T(\mathcal{X})$ by
$\mu = \int_{\mathscr{M}^T(\mathcal{X})} \nu \, d\lambda(\nu)$.
For any positive number $\varepsilon$, we have 
\begin{equation*}
   \begin{split}
    &   R(\mathbf{d}, \mu, \varepsilon) \\
   &       =  \inf\left\{\int_{\mathscr{M}^T(\mathcal{X})} R(\mathbf{d},\nu, \varepsilon_\nu)d\lambda(\nu) \middle|\, 
\text{\parbox{8cm}{$\mathscr{M}^T(\mathcal{X})\ni \nu \mapsto \varepsilon_\nu \in (0, \infty)$ is a measurable function with  
$\int_{\mathscr{M}^T(\mathcal{X})} \varepsilon_\nu \, d\lambda(\nu) \leq \varepsilon$}}\right\}. 
   \end{split}
\end{equation*}   
\end{theorem}

\begin{remark}
So far we have considered only $\mathbb{R}^d$-actions.
However we notice that the statement of Theorem \ref{theorem: rate distortion function and convex combination of measures}
holds for $\mathbb{Z}^d$-actions as well.
This might be useful in a future study of geometric applications of rate distortion theory.
\end{remark}

\subsection{Proofs of Theorems \ref{theorem: upper rate distortion dimension and convex combination} and 
\ref{theorem: lower rate distortion dimension and convex combination}}
\label{subsection: proofs of main theorems}

Let $(\mathcal{X}, \mathbf{d})$ be a compact metric space with a continuous action 
$T\colon \mathbb{R}^d\times \mathcal{X}\to \mathcal{X}$.
Let $\lambda$ be a Borel probability measure on $\mathscr{M}^T(\mathcal{X})$ and define a measure 
$\mu \in \mathscr{M}^T(\mathcal{X})$ by
$\mu = \int_{\mathscr{M}^T(\mathcal{X})} \nu \, d\lambda(\nu)$.
We prove Theorems \ref{theorem: upper rate distortion dimension and convex combination} and 
\ref{theorem: lower rate distortion dimension and convex combination} below.

\begin{theorem}[$=$ Theorem \ref{theorem: upper rate distortion dimension and convex combination}]
If the upper metric mean dimension $\umdimm(\mathcal{X}, T, \mathbf{d})$ is finite then 
\[ \urdim(\mathcal{X}, T, \mathbf{d}, \mu) \leq 
    \int_{\mathscr{M}^T(\mathcal{X})} \urdim(\mathcal{X}, T, \mathbf{d}, \nu) \, d\lambda(\nu). \]
\end{theorem}

\begin{proof}
Recall that 
\[  \umdimm(\mathcal{X}, T, \mathbf{d}) = \limsup_{\varepsilon \to 0} \frac{S(\mathcal{X}, T, \mathbf{d}, \varepsilon)}{\log (1/\varepsilon)}. \]

By Theorem \ref{theorem: rate distortion function and convex combination of measures}, for any $\varepsilon>0$ we have 
\[  R(\mathbf{d},\mu, \varepsilon) \leq \int_{\mathscr{M}^T(\mathcal{X})} R(\mathbf{d}, \nu, \varepsilon)\, d\lambda(\nu). \]
Hence 
\[  \frac{R(\mathbf{d},\mu, \varepsilon)}{\log (1/\varepsilon)} \leq 
      \int_{\mathscr{M}^T(\mathcal{X})} \frac{R(\mathbf{d}, \nu, \varepsilon)}{\log (1/\varepsilon)} \, d\lambda(\nu). \]
From (\ref{eq: rate distortion function and varepsilon entropy}) in \S \ref{subsection: rate distortion function},
\[  \frac{R(\mathbf{d}, \nu, \varepsilon)}{\log (1/\varepsilon)}  \leq  \frac{S(\mathcal{X}, T, \mathbf{d}, \varepsilon)}{\log (1/\varepsilon)}. \]
From $\umdimm(\mathcal{X}, T, \mathbf{d}) < \infty$, we have 
\[  \sup_{\varepsilon>0}  \frac{R(\mathbf{d}, \nu, \varepsilon)}{\log (1/\varepsilon)}  \leq
    \sup_{\varepsilon>0} \frac{S(\mathcal{X}, T, \mathbf{d},\varepsilon)}{\log (1/\varepsilon)} < \infty. \]
Namely, $\frac{R(\mathbf{d}, \nu, \varepsilon)}{\log (1/\varepsilon)}$ is a bounded function of $\nu$ and $\varepsilon$.
Then we can use Fatou’s lemma and conclude
\[  \limsup_{\varepsilon \to 0} \frac{R(\mathbf{d},\mu, \varepsilon)}{\log (1/\varepsilon)}  \leq 
   \limsup_{\varepsilon\to 0}  \int_{\mathscr{M}^T(\mathcal{X})} \frac{R(\mathbf{d}, \nu, \varepsilon)}{\log (1/\varepsilon)} \, d\lambda(\nu)
   \leq  \int_{\mathscr{M}^T(\mathcal{X})}
    \limsup_{\varepsilon\to 0} \frac{R(\mathbf{d}, \nu, \varepsilon)}{\log (1/\varepsilon)} \, d\lambda(\nu). \]
Thus we have 
$\urdim(\mathcal{X}, T, \mathbf{d}, \mu) \leq 
    \int_{\mathscr{M}^T(\mathcal{X})} \urdim(\mathcal{X}, T, \mathbf{d}, \nu) \, d\lambda(\nu)$.
\end{proof}

\begin{theorem}[$=$ Theorem \ref{theorem: lower rate distortion dimension and convex combination}]
We have
\[  \lrdim(\mathcal{X}, T, \mathbf{d}, \mu) \geq 
    \int_{\mathscr{M}^T(\mathcal{X})} \lrdim(\mathcal{X}, T, \mathbf{d}, \nu) \, d\lambda(\nu). \]
\end{theorem}

\begin{proof}
We take a decreasing sequence of positive numbers $\varepsilon_1>\varepsilon_2>\varepsilon_3>\dots$ such that 
\[  \varepsilon_n \leq 2^{-n}, \quad \lrdim(\mathcal{X}, T, \mathbf{d}, \mu) = \lim_{n\to \infty} 
     \frac{R(\mathbf{d}, \mu, \varepsilon_n)}{\log (1/\varepsilon_n)}. \]
By Theorem \ref{theorem: rate distortion function and convex combination of measures}, for each $n\geq 1$, we can find 
a measurable function 
$\mathscr{M}^T(\mathcal{X})\ni \nu \mapsto \varepsilon_{n,\nu} \in (0, \infty)$ such that 
$\int_{\mathscr{M}^T(\mathcal{X})} \varepsilon_{n,\nu}\, d\lambda(\nu) \leq \varepsilon_n$ and 
\begin{equation} \label{eq: choice of varepsilon_nu(n)}
    R(\mathbf{d}, \mu, \varepsilon_n) \geq \int_{\mathscr{M}^T(\mathcal{X})} R\left(\mathbf{d}, \nu, \varepsilon_{n,\nu}\right) \, d\lambda(\nu) -1. 
\end{equation}    

Take $0<c<1$ and let 
$A_n = \{\nu \in \mathscr{M}^T(\mathcal{X})\mid \varepsilon_{n,\nu} < \varepsilon_n^c\}$.
We have $\lambda\left(\mathscr{M}^T(\mathcal{X})\setminus A_n\right) \leq \varepsilon_n^{1-c} \leq 2^{-n(1-c)}$.
Hence the sum of $\lambda\left(\mathscr{M}^T(\mathcal{X})\setminus A_n\right)$ over $n\geq 1$ converges.
By the first Borel--Canteli lemma, 
\[  \lambda \left(\bigcap_{n=1}^\infty \bigcup_{k\geq n} \left(\mathscr{M}^T(\mathcal{X})\setminus A_k\right)\right) = 0. \]
Set $B = \bigcup_{n=1}^\infty \bigcap_{k\geq n} A_k$. We have $\lambda(B) = 1$.
For any $\nu \in B$, we have $\varepsilon_{n, \nu} < \varepsilon_n^c$ for all but finitely many $n$.

By (\ref{eq: choice of varepsilon_nu(n)})
 \begin{align*}
    \frac{R(\mathbf{d}, \mu, \varepsilon_n)}{\log (1/\varepsilon_n)} &\geq \int_{\mathscr{M}^T(\mathcal{X})} 
    \frac{R(\mathbf{d}, \nu, \varepsilon_{n,\nu})}{\log (1/\varepsilon_n)} d\lambda(\nu) - \frac{1}{\log (1/\varepsilon_n)} \\
    & =  \int_{B} \frac{R(\mathbf{d}, \nu, \varepsilon_{n,\nu})}{\log (1/\varepsilon_n)}d\lambda(\nu) 
    - \frac{1}{\log (1/\varepsilon_n)} \quad 
    (\text{by $\lambda(B)=1$}). 
  \end{align*}
We notice that rate distortion function is always non-negative.
Hence we can use Fatou’s lemma and get 
\[ \lrdim(\mathcal{X}, T, \mathbf{d}, \mu) = \lim_{n\to \infty} 
     \frac{R(\mathbf{d}, \mu, \varepsilon_n)}{\log (1/\varepsilon_n)}
     \geq \int_{B} \liminf_{n\to \infty} \frac{R(\mathbf{d}, \nu, \varepsilon_{n,\nu})}{\log (1/\varepsilon_n)} \, d\lambda(\nu). \]
For $\nu\in B$ we have 
\[ \frac{R(\mathbf{d}, \nu, \varepsilon_{n,\nu})}{\log (1/\varepsilon_n)}  
    \geq  \frac{R(\mathbf{d}, \nu, \varepsilon_n^c)}{\log (1/\varepsilon_n)} 
     = c\frac{R(\mathbf{d}, \nu, \varepsilon_n^c)}{\log (1/\varepsilon_n^c)} \quad \text{for all but finitely many $n\geq 1$}. \]
Therefore 
  \begin{align*}
    \int_{B} \liminf_{n\to \infty} \frac{R(\mathbf{d}, \mu, \varepsilon_{n,\nu})}{\log (1/\varepsilon_n)} \, d\lambda(\nu) & \geq 
    c \int_B \liminf_{n\to \infty} \frac{R(\mathbf{d}, \nu, \varepsilon_n^c)}{\log (1/\varepsilon_n^c)} d\lambda(\nu) \\
    & \geq  c \int_B \lrdim(\mathcal{X}, T, \mathbf{d}, \nu) \, d\lambda(\nu) \\
   &  = c \int_{\mathscr{M}^T(\mathcal{X})} \lrdim(\mathcal{X}, T, \mathbf{d}, \nu) \, d\lambda(\nu).  
  \end{align*}
In the last line we have used $\lambda(B)=1$. Thus we have 
\[  \lrdim(\mathcal{X}, T, \mathbf{d}, \mu) \geq  c \int_{\mathscr{M}^T(\mathcal{X})} \lrdim(\mathcal{X}, T, \mathbf{d}, \nu) \, d\lambda(\nu).   \]
Letting $c\to 1$, we get the conclusion.
\end{proof}

\section{Construction of an example: Proof of Theorem \ref{theorem: irregular behavior of upper rate distortion dimension}}
\label{section: construction of an example}

The purpose of this section is to prove Theorem \ref{theorem: irregular behavior of upper rate distortion dimension}.
Let $(V, \norm{\cdot})$ be an infinite dimensional Banach space.
For example, $V=\ell^2(\mathbb{N})$ or $\ell^\infty(\mathbb{N})$ will work well.
Let $B = \{v\in V\mid \norm{v}\leq 1\}$ be the unit ball of $V$.
Let $C(\mathbb{R}^d, B)$ be the space of all continuous maps 
$x\colon \mathbb{R}^d\to B$.
We define a metric $\mathbf{d}$ on $C(\mathbb{R}^d, B)$ by 
\[ \mathbf{d}(x, y) = \sup_{n\geq 1} \left(2^{-n} \sup_{|t|_\infty \leq n} \norm{x(t)-y(t)}\right). \]
Here $|t|_\infty := \max(|t_1|, |t_2|, \dots, |t_d|)$ denotes the max-norm of $t=(t_1, \dots, t_d)\in \mathbb{R}^d$.
The group $\mathbb{R}^d$ continuously acts on $C(\mathbb{R}^d, B)$ by 
\[  T\colon \mathbb{R}^d\times C(\mathbb{R}^d, B) \to C(\mathbb{R}^d, B), \quad 
     \left(s, x(\cdot)\right) \mapsto x(\cdot+s). \]
Namely $(T^sx)(t) = x(t+s)$.

\begin{figure}[h] 
    \centering
    \includegraphics[width=3.0in]{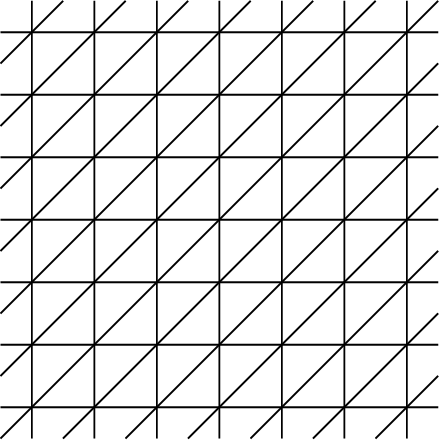}
    \caption{A $\mathbb{Z}^2$-invariant triangulation of the plane. 
    The vertexes are points of the standard lattice $\mathbb{Z}^2$.}    \label{figure: triangulation}
\end{figure}

We take a triangulation (simplicial complex structure) $\Gamma$ of $\mathbb{R}^d$ such that 
\begin{itemize}
  \item the vertex sets of $\Gamma$ is the standard lattice $\mathbb{Z}^d$,
  \item $\Gamma$ is invariant under the natural action of $\mathbb{Z}^d$ on $\mathbb{R}^d$, namely 
  $u+\Gamma = \Gamma$ for all $u\in \mathbb{Z}^d$.
\end{itemize}
Figure \ref{figure: triangulation} shows an example of such $\Gamma$ for the case of $d=2$.

A continuous map $x\colon \mathbb{R}^d\to B$ is said to be \textbf{$\Gamma$-piecewisely linear} 
if, for every simplex $\Delta$ of $\Gamma$, we have
\[ x\left(\sum_{i=0}^k t_i e_i\right) = \sum_{i=0}^k t_i x(e_i)  \quad \text{for $t_i\geq 0$ with $\sum_{i=0}^k t_i = 1$},  \]
where $e_0, \dots, e_k$ are the vertices of $\Delta$.
A $\Gamma$-picewisely linear map $x$ is uniquely determined by its values on the vertex set $\mathbb{Z}^d$.

\begin{lemma} \label{lemma: Lipschitz constant of PL maps}
\begin{enumerate}
\item There is $C>0$ such that every $\Gamma$-piecewisely linear map $x\colon \mathbb{R}^d\to B$ is $C$-Lipschitz:
\[  \norm{x(t)- x(u)} \leq C |t-u|_\infty. \]
\item For $L>0$, $s, t\in \mathbb{R}^d$ and $\Gamma$-piecewisely linear map $x\colon \mathbb{R}^d\to B$ we have
\[  \mathbf{d}_L(T^s x, T^t x)\leq \frac{C}{2} |s-t|_\infty. \]
\end{enumerate}
\end{lemma}

\begin{proof}
(1)
Since the triangulation $\Gamma$ is $\mathbb{Z}^d$-invariant, we can find $C>0$ such that, for every simplex $\Delta$ of 
$\Gamma$ of vertices $e_0, \dots, e_k$, we have 
\[  \sum_{i=0}^k |t_i-u_i|  \leq  C \left|\sum_{i=0}^k t_i e_i - \sum_{i=0}^k u_i e_i\right|_\infty \]
where $t_i, u_i$ are any non-negative numbers with $\sum_{i=0}^k t_i = \sum_{i=0}^k u_i = 1$.
Then for any two points $t=\sum_{i=0}^k t_i e_i$ and $u=\sum_{i=0}^k u_i e_i$ of $\Delta$ we have 
\[  \norm{x(t)-x(u)} \leq  \sum_{i=0}^k |t_i-u_i|  \leq  C \left|\sum_{i=0}^k t_i e_i - \sum_{i=0}^k u_i e_i\right|_\infty = C|t-u|_\infty, \]
where we have used $\norm{x(e_i)}\leq 1$ in the first inequality.

In general, let $t$ and $u$ be any two points of $\mathbb{R}^d$. We consider the line segment 
$p(s) := (1-s)t + su$ $(0\leq s \leq 1)$ between $t$ and $u$.
We can find $0=s_0 < s_1 < s_2<\dots < s_n=1$ such that each pair $p(s_j), p(s_{j+1})$ $(0\leq j \leq n-1)$ 
belong to a common simplex of $\Gamma$.
Then 
\begin{align*}
  \norm{x(t)-x(u)} &\leq \sum_{j=0}^{n-1} \norm{x\left(p(s_j)\right) - x\left(p(s_{j+1}\right)}  \\
  &\leq C \sum_{j=0}^{n-1} \left|p(s_j) - p(s_{j+1})\right|_\infty \\
  & = C \sum_{j=0}^{n-1} (s_{j+1}-s_j)|t-u|_\infty \\
  & = C |t-u|_\infty.
\end{align*}  
(2) is a direct consequence of (1) and the definitions.
\end{proof}

For natural numbers $m\geq 1$ we take finite subsets $A_m\subset B$ satisfying the following conditions.
\begin{enumerate}
  \item $\norm{v-w}=\frac{1}{m}$ for all distinct $v, w\in A_m$.
  \item $0\not \in A_m$ and $\max_{v \in A_m} \norm{a} \to 0$ as $m\to \infty$.  
  \item $A_m\cap A_{m^\prime} =\emptyset$ for $m\neq m^\prime$, 
           and $\bigcup_{m=1}^\infty A_m$ is a linearly independent subset of $V$.  
  \item $|A_m| = 2^{3^m}$, where $|A_m|$ denotes the cardinality of $A_m$.
\end{enumerate}
We can find such $A_m$ if $V=\ell^2(\mathbb{N})$ or $\ell^\infty(\mathbb{N})$.

We define a subset $\mathcal{Y}_m \subset C(\mathbb{R}^d, B)$ as the set of all 
$\Gamma$-piecewisely linear maps $x\colon \mathbb{R}^d\to B$ satisfying 
$x(\mathbf{n}) \in A_m$  for all $\mathbf{n} \in \mathbb{Z}^d$.
This is a compact subset and invariant under $T|_{\mathbb{Z}^d}$ (the restriction of $T$ to $\mathbb{Z}^d$).
We define $\mathcal{X}_m \subset C(\mathbb{R}^d, B)$ by 
\[  \mathcal{X}_m = \bigcup_{s\in [0, 1)^d} T^{s} \mathcal{Y}_m. \]
This is a compact $T$-invariant subset of $C(\mathbb{R}^d, B)$.
We also set $\mathcal{X}_0 = \{0\}$ where $0$ denotes the zero function in $C(\mathbb{R}^d, B)$.
We have $\mathcal{X}_m\cap \mathcal{X}_{m^\prime} = \emptyset$ for $m\neq m^\prime$ by the condition (3) above.
Finally we define $\mathcal{X}\subset C(\mathbb{R}^d, B)$ by 
\[  \mathcal{X} = \bigcup_{m=0}^\infty \mathcal{X}_m. \]
$\mathcal{X}$ is a compact $T$-invariant subset of $C(\mathbb{R}^d, B)$.
Here the compactness follows from the condition that $\max_{v \in A_m} \norm{v} \to 0$ $(m\to \infty)$.
We will show that $(\mathcal{X}, T, \mathbf{d})$ satisfies the statement of 
Theorem \ref{theorem: irregular behavior of upper rate distortion dimension}.

\begin{lemma} \label{lemma: topological entropy of X_m}
For every $m\geq 0$, the topological entropy of $(\mathcal{X}_m, T)$ is finite and hence
we have $\umdimm\left(\mathcal{X}_m, T, \mathbf{d}\right) = 0$.
\end{lemma}

\begin{proof}
For $m=0$, this is obvious because $\mathcal{X}_0 = \{0\}$ is a single point.
So we assume $m\geq 1$.
Let $\varepsilon>0$ and $L>0$.
From Lemma \ref{lemma: Lipschitz constant of PL maps} (2)
  \begin{align*}
    \#\left(\mathcal{X}_m, \mathbf{d}_L, \varepsilon\right) & \leq 
    \sum_{s\in [0, 1)^d\cap \left(\frac{\varepsilon}{2C}\mathbb{Z}^d\right)}
    \#\left(T^s \mathcal{Y}_m, \mathbf{d}_{L}, \frac{\varepsilon}{2}\right)  \\
    & \leq   \sum_{s\in [0, 1)^d\cap \left(\frac{\varepsilon}{2C}\mathbb{Z}^d\right)}
   \#\left(\mathcal{Y}_m, \mathbf{d}_{L+1}, \frac{\varepsilon}{2}\right) \\
    & \leq  \left(1+\frac{2C}{\varepsilon}\right)^d   \#\left(\mathcal{Y}_m, \mathbf{d}_{L+1}, \frac{\varepsilon}{2}\right).
  \end{align*}
Choose $\ell=\ell(\varepsilon)>0$ satisfying $2^{-\ell} < \varepsilon/2$. Then
\[ \#\left(\mathcal{Y}_m, \mathbf{d}_{L+1}, \frac{\varepsilon}{2}\right) \leq  \left|A_m\right|^{(L+2\ell+1)^d}. \]
Hence
\[ S(\mathcal{X}_m, T, \mathbf{d}, \varepsilon) 
   = \lim_{L\to \infty} \frac{\log \#\left(\mathcal{X}_m, \mathbf{d}_L, \varepsilon\right)}{L^d}
   \leq \log |A_m|. \]
Thus 
\[ h_{\mathrm{top}}\left(\mathcal{X}_m, T\right) = \lim_{\varepsilon \to 0} S(\mathcal{X}_m, T, \mathbf{d}, \varepsilon) \leq 
    \log |A_m| < \infty,  \]
\[  \umdimm\left(\mathcal{X}_m, T, \mathbf{d}\right) = \limsup_{\varepsilon \to 0} 
      \frac{S(\mathcal{X}_m, T, \mathbf{d}, \varepsilon)}{\log (1/\varepsilon)} = 0. \]
\end{proof}

\begin{lemma} \label{lemma: rate distortion dimension of ergodic measures}
We have $\rdim\left(\mathcal{X},T, \mathbf{d}, \nu\right) = 0$
for every ergodic probability measure $\nu$ of $(\mathcal{X}, T)$ .
\end{lemma}

\begin{proof}
We have a decomposition $\mathcal{X} = \bigcup_{m=0}^\infty \mathcal{X}_m$.
The sets $\mathcal{X}_m$ are mutually disjoint, compact and $T$-invariant subsets.
We have 
\[ 1 = \nu(\mathcal{X}) = \sum_{m=0}^\infty \nu(\mathcal{X}_m). \]
Then $\nu(\mathcal{X}_m) = 1$ for some $m\geq 0$ because $\nu$ is ergodic.
From Lemma \ref{lemma: topological entropy of X_m}
\[ \rdim\left(\mathcal{X},T, \mathbf{d}, \nu\right) = \rdim\left(\mathcal{X}_m, T, \mathbf{d}, \nu\right)  \leq  
\umdimm\left(\mathcal{X}_m, T, \mathbf{d}\right) =0. \]
\end{proof}

For $m\geq 1$, let $\Phi_m\colon A_m^{\mathbb{Z}^d}\to  \mathcal{Y}_m$ be the “natural” map.
Namely, for $\mathbf{a} = (a_{\mathbf{n}})_{\mathbf{n}\in \mathbb{Z}^d} \in A_m^{\mathbb{Z}^d}$, we define 
$\Phi_m(\mathbf{a})$ to be the (unique) $\Gamma$-piecewisely linear map $x\colon \mathbb{R}^d\to B$ satisfying 
$x(\mathbf{n}) = a_{\mathbf{n}}$ for all $\mathbf{n}\in \mathbb{Z}^d$.
Let $p_m$ be the uniform probability measure on $A_m$.
Let $p_m^{\otimes \mathbb{Z}^d}$ be the product measure on $A_m^{\mathbb{Z}^d}$ and set 
\[ \nu_m = (\Phi_m)_* \left(p_m^{\otimes \mathbb{Z}^d}\right). \]
This is a $T|_{\mathbb{Z}^d}$-invariant probability measure on $\mathcal{Y}_m$.
We define a $T$-invariant probability measure $\mu_m$ on $\mathcal{X}_m$ by 
\[ \mu_m = \int_{[0, 1)^d} T^s_*\nu_m\, d\mathbf{m}(s). \]
Finally we define $\mu \in \mathscr{M}^T(\mathcal{X})$ by 
\[ \mu = \sum_{m=1}^\infty 2^{-m} \mu_m. \]

We will prove that $\rdim\left(\mathcal{X}, T, \mathbf{d}, \mu\right) = \infty$.
The strategy is to relate $R(\mathbf{d},\mu,\varepsilon)$ to a convex combination of 
$R(\mathbf{d},\mu_n,\varepsilon)$ (Theorem \ref{theorem: rate distortion function and convex combination of measures})
and establish an appropriate lower bound on $R(\mathbf{d}, \mu_m, \varepsilon)$ by using 
Proposition \ref{proposition: lower bound on rate distortion function}.
For applying that proposition, we need to estimate the integral 
\[  \int_{\mathcal{X}} 2^{-\frac{a}{L^d}\int_{[0, L)^d}\mathbf{d}(T^t x, y_t)d\mathbf{m}(t)}d\mu_m(x) 
     =  \int_{[0, 1)^d} \left(\int_{T^s\mathcal{Y}_m} 2^{-\frac{a}{L^d}\int_{[0, L)^d}\mathbf{d}(T^t x, y_t)d\mathbf{m}(t)}
         d\left(T^s_*\nu_m(x)\right)\right) ds \]
for $L>0$, $y\in L^1\left([0, L)^d, \mathcal{X}\right)$ and an appropriately chosen positive number $a$.
(We will see below that the choice $a:=4m L^d 3^m$ will work.)
Here is a remark on the notation: For $y\in L^1([0, L)^d, \mathcal{X})$ and $t\in [0, L)^d$, we denote the value of $y$ at $t$ by 
$y_t$. So $y_t$ is a continuous map from $\mathbb{R}^d$ to $B$.
For $u\in \mathbb{R}^d$, we denote the value of $y_t$ at the point $u$ by $y_t(u)$.

\begin{lemma} \label{lemma: distortion estimate}
Let $L\in \mathbb{N}$, $s\in [0, 1)^d$, $x\in \mathcal{X}$ and $y\in L^1\left([0, L)^d, \mathcal{X}\right)$.
We have 
\[  \int_{[0, L)^d} \mathbf{d}\left(T^t x, y_t\right) d \mathbf{m}(t) \geq 
     \frac{1}{2}\sum_{\mathbf{n}\in \mathbb{Z}^d\cap [0, L)^d}
     \norm{x(s+\mathbf{n}) - \int_{\mathbf{n}+[0, 1)^d} y_t(s+\mathbf{n}-t)\, d\mathbf{m}(t)}. \]
\end{lemma}

\begin{proof}
We have 
\[ \mathbf{d}(T^t x, y_t) = \sup_{n\geq 1}\left( 2^{-n}\sup_{|u|_\infty \leq n} \norm{T^t x(u) - y_t(u)} \right)
                                     \geq \frac{1}{2} \max_{|u|_\infty \leq 1} \norm{x(u+t)- y_t(u)}. \]
For $\mathbf{n}\in [0, L)^d\cap \mathbb{Z}^d$ and $t\in \mathbf{n}+[0, 1)^d$, letting $u=s+\mathbf{n}-t\in [-1,1]^d$, we get 
\[  \mathbf{d}(T^t x, y_t) \geq \frac{1}{2} \norm{x(s+\mathbf{n})- y_t(s+\mathbf{n}-t)}. \]
Integrating this over $t\in \mathbf{n}+[0, 1)^d$,
\begin{align*}
    \int_{\mathbf{n}+[0, 1]^d} \mathbf{d}(T^t x, y_t) d\mathbf{m}(t)  & \geq 
    \frac{1}{2} \int_{\mathbf{n}+[0,1)^d} \norm{x(s+\mathbf{n})- y_t(s+\mathbf{n}-t)} d\mathbf{m}(t) \\
    & \geq \frac{1}{2} \norm{x(s+\mathbf{n}) - \int_{\mathbf{n}+[0,1)^d} y_t(s+\mathbf{n}-t) \, d\mathbf{m}(t)}.
\end{align*}    
Summing this over $\mathbf{n}\in [0, L)^d\cap \mathbb{Z}^d$, we get the claim of the lemma.
\end{proof}

Let  $L\in \mathbb{N}$, $s\in [0, 1)^d$ and $y\in L^1\left([0, L)^d, \mathcal{X}\right)$.
We set 
\[ z_s(\mathbf{n})  =  \int_{\mathbf{n}+[0,1)^d} y_t(s+\mathbf{n}-t)\, d\mathbf{m}(t) \in B \quad 
    (\mathbf{n}\in \mathbb{Z}^d\cap [0, L)^d). \]
From Lemma \ref{lemma: distortion estimate}, for any positive number $a$
\begin{align*}
    \int_{T^s\mathcal{Y}_m} 2^{-\frac{a}{L^d}\int_{[0, L)^d}\mathbf{d}(T^t x, y_t)d\mathbf{m}(t)}
         d\left(T^s_*\nu_m(x)\right) & \leq 
     \int_{T^s \mathcal{Y}_m}  2^{-\frac{a}{2L^d}\sum_{\mathbf{n}\in [0, L)^d\cap \mathbb{Z}^d}\norm{x(s+\mathbf{n})-z_s(\mathbf{n})}}
     d\left(T^s_*\nu_m(x)\right)  \\
     & = \int_{T^s \mathcal{Y}_m} \prod_{\mathbf{n}\in [0, L)^d\cap \mathbb{Z}^d}
           2^{-\frac{a}{2L^d}\norm{x(s+\mathbf{n}) - z_s(\mathbf{n})}} d\left(T^s_*\nu_m(x)\right).
\end{align*}           
When $x\in T^s \mathcal{Y}_m$ is distributed according to $T^s_*\nu_m$, the point 
$\left(x(s+\mathbf{n})\right)_{\mathbf{n}\in [0, L)^d\cap \mathbb{Z}^d}$ is distributed according to the product measure 
$\prod_{\mathbf{n}\in [0, L)^d\cap \mathbb{Z}^d} p_m$.
(Recall that $p_m$ is the uniform measure on the finite set $A_m$.)
Hence
\[  \int_{T^s \mathcal{Y}_m} \prod_{\mathbf{n}\in [0, L)^d\cap \mathbb{Z}^d}
           2^{-\frac{a}{2L^d}\norm{x(s+\mathbf{n}) - z_s(\mathbf{n})}} d\left(T^s_*\nu_m(x)\right)
     =     \prod_{\mathbf{n}\in [0, L)^d\cap \mathbb{Z}^d} \int_{A_m} 2^{-\frac{a}{2L^d}\norm{v-z_s(\mathbf{n})}} d p_m(v). \]   
Therefore 
\[ \int_{T^s\mathcal{Y}_m} 2^{-\frac{a}{L^d}\int_{[0, L)^d}\mathbf{d}(T^t x, y_t)d\mathbf{m}(t)}
         d\left(T^s_*\nu_m(x)\right)  \leq 
 \prod_{\mathbf{n}\in [0, L)^d\cap \mathbb{Z}^d} \int_{A_m} 2^{-\frac{a}{2L^d}\norm{v-z_s(\mathbf{n})}} d p_m(v). \]  

Recall that we assumed $\norm{v-w}=\frac{1}{m}$ for all distinct $v, w\in A_m$.
It follows that $\norm{v-z_s(\mathbf{n})}\geq \frac{1}{2m}$ for all but one $v\in A_m$.
Therefore 
\[  \int_{A_m} 2^{-\frac{a}{2L^d}\norm{v-z_s(\mathbf{n})}} d p_m(v) \leq 
     \frac{1}{|A_m|}\left\{1+(|A_m|-1)2^{-\frac{a}{4mL^d}}\right\} \leq 
      \frac{1}{|A_m|} + 2^{-\frac{a}{4mL^d}}. \]
We recall that $|A_m|=2^{3^m}$. 
Then the right-most hand is equal to $2^{-3^m} + 2^{-\frac{a}{4mL^d}}$.
Hence
\[ \int_{T^s\mathcal{Y}_m} 2^{-\frac{a}{L^d}\int_{[0, L)^d}\mathbf{d}(T^t x, y_t)d\mathbf{m}(t)}
         d\left(T^s_*\nu_m(x)\right)  \leq \left(2^{-3^m} + 2^{-\frac{a}{4mL^d}}\right)^{L^d}. \]
Integrating this over $s\in [0, 1)^d$, we have 
\[ \int_{\mathcal{X}} 2^{-\frac{a}{L^d}\int_{[0, L)^d}\mathbf{d}(T^t x, y_t)d\mathbf{m}(t)}d\mu_m(x) 
    \leq  \left(2^{-3^m} + 2^{-\frac{a}{4mL^d}}\right)^{L^d}. \]
Now we choose $a:=4m L^d 3^m$. 
Then the right-hand side is $(2^{-3^m} + 2^{-3^m})^{L^d} = 2^{(1-3^m)L^d}$.
Namely 
\[ \int_{\mathcal{X}} 2^{-4m\cdot 3^m \int_{[0, L)^d}\mathbf{d}(T^t x, y_t)d\mathbf{m}(t)}d\mu_m(x) \leq 2^{(1-3^m)L^d}. \]
In other words, for $a= 4m L^d 3^m$ and $\lambda := 2^{(3^m-1)L^d}$, we have 
\[ \int_{\mathcal{X}}\lambda 2^{-\frac{a}{L^d}\int_{[0, L)^d}\mathbf{d}(T^t x, y_t)d\mathbf{m}(t)}d\mu_m(x)  \leq 1  \]
for all $y\in L^1\left([0, L)^d, \mathcal{X}\right)$.
Then we can use Proposition \ref{proposition: lower bound on rate distortion function} and get 
\[  R_L(\mathbf{d}, \mu_m, \varepsilon) \geq  -a \varepsilon + \log \lambda  = -4m\varepsilon L^d 3^m + (3^m-1) L^d \quad 
     (\varepsilon >0). \]
Dividing this by $L^d$ and letting $L\to \infty$,
\[  R(\mathbf{d}, \mu_m, \varepsilon) \geq  -4m\varepsilon 3^m + 3^m -1. \]     
For $0<\varepsilon < 3^{-m}$, the right-hand side is greater than $3^m - 4m -1$.
Therefore we conclude 

\begin{proposition} \label{proposition: lower bound on rate distortion function of mu_m}
For $0<\varepsilon<3^{-m}$ we have $R(\mathbf{d}, \mu_m, \varepsilon) \geq 3^m - 4m -1$.
\end{proposition}

Now we apply Theorem \ref{theorem: rate distortion function and convex combination of measures} to 
$\mu = \sum_{m=1}^\infty 2^{-m}\mu_m$: 
\[ R(\mathbf{d},\mu, \varepsilon) = \inf\left\{\sum_{m=1}^\infty 2^{-m} R(\mathbf{d}, \mu_m, \varepsilon_m)\middle|\, 
     \varepsilon_m>0 \, (m\in \mathbb{N}) \text{ with } \sum_{m=1}^\infty 2^{-m}\varepsilon_m \leq \varepsilon\right\}. \]
Let $n$ be a natural number and assume $0<\varepsilon < 6^{-n}$.
If a sequence of positive numbers $\varepsilon_m$ satisfies $\sum_{m=1}^\infty 2^{-m}\varepsilon_m \leq \varepsilon$, then 
$2^{-n}\varepsilon_n \leq \varepsilon < 6^{-n}$ and hence $\varepsilon_n < 3^{-n}$.
By Proposition \ref{proposition: lower bound on rate distortion function of mu_m}
\[  2^{-n}  R(\mathbf{d}, \mu_n, \varepsilon_n)  \geq  \left(\frac{3}{2}\right)^n - n 2^{2-n} - 2^{-n}. \]
Therefore 
\[  R(\mathbf{d},\mu,\varepsilon) \geq   \left(\frac{3}{2}\right)^n - n 2^{2-n} - 2^{-n} \quad 
     (0<\varepsilon < 6^{-n}). \]

\begin{proposition} \label{proposition: rate distortion dimension is infinite}
 $\rdim\left(\mathcal{X}, T, \mathbf{d}, \mu\right) = \infty$.
\end{proposition}

\begin{proof}
Let $0<\varepsilon<1/6$. We choose a natural number $n$ with $6^{-n-1}\leq \varepsilon < 6^{-n}$.
Then we have 
\[ R(\mathbf{d}, \mu, \varepsilon) \geq  \left(\frac{3}{2}\right)^n - n 2^{2-n} - 2^{-n}, \quad 
    \log (1/\varepsilon) \leq  (n+1)\log 6. \]
Hence 
\[  \frac{R(\mathbf{d}, \mu, \varepsilon)}{\log (1/\varepsilon)} \geq \frac{\left(\frac{3}{2}\right)^n - n 2^{2-n} - 2^{-n}}{(n+1)\log 6}. \]
We have $n\to \infty$ as $\varepsilon \to 0$. Therefore 
\[   \rdim\left(\mathcal{X}, T, \mathbf{d}, \mu\right) 
      = \lim_{\varepsilon \to 0} \frac{R(\mathbf{d}, \mu, \varepsilon)}{\log (1/\varepsilon)} = \infty. \]
\end{proof}

By combining Lemmas \ref{lemma: rate distortion dimension of ergodic measures} and 
\ref{proposition: rate distortion dimension is infinite}, we have shown that $(\mathcal{X}, T, \mathbf{d})$ satisfies 
\begin{itemize} 
   \item every ergodic measure $\nu$ on $(\mathcal{X}, T)$ has zero rate distortion dimension, 
   \item there exists an invariant probability measure $\mu$ on $(\mathcal{X}, T)$ of infinite rate distortion dimension.
\end{itemize}   
This proves Theorem \ref{theorem: irregular behavior of upper rate distortion dimension}.

\begin{remark}
By modifying the above construction, we can also prove the following statement:
Let $c$ be an arbitrary nonnegative real number (including $c=\infty$).
There exists a compact metric space $(\mathcal{X}, \mathbf{d})$ with a continuous action 
$T\colon \mathbb{R}^d\times \mathcal{X}\to \mathcal{X}$ satisfying the following two conditions.
  \begin{enumerate}
     \item Every ergodic measure $\nu\in \mathscr{M}^T_{\mathrm{erg}}(\mathcal{X})$ satisfies 
              $\rdim\left(\mathcal{X}, T, \mathbf{d}, \nu\right) = 0$.
     \item There exists $\mu \in \mathscr{M}^T(\mathcal{X})$ satisfying $\rdim\left(\mathcal{X}, T, \mathbf{d}, \mu\right) = c$.
  \end{enumerate}
\end{remark}

\section{Another example}  \label{section: another example}

One of our main theorems, Theorem \ref{theorem: upper rate distortion dimension and convex combination}, claims that 
if a continuous action $T\colon \mathbb{R}^d\times \mathcal{X}\to \mathcal{X}$ on a compact metric space $(\mathcal{X}, \mathbf{d})$ 
has a finite upper metric mean dimension and if $\lambda$ is a Borel probability measure on $\mathscr{M}^T(\mathcal{X})$, then 
we have
\begin{equation} \label{eq: convexity inequality for upper rate distortion dimension}
    \urdim \left(\mathcal{X}, T, \mathbf{d}, \int_{\mathscr{M}^T(\mathcal{X})}\nu \, d\lambda(\nu)\right) 
    \leq  \int_{\mathscr{M}^T(\mathcal{X})} \urdim\left(\mathcal{X}, T, \mathbf{d}, \nu\right) d\lambda(\nu). 
\end{equation}    
The purpose of this section is to show that this inequality may be a strict inequality in general.
Indeed we will prove the following statement.

\begin{proposition} \label{proposition: convexity inequality may be strict}
   There exists a continuous action $T\colon \mathbb{R}^d\times \mathcal{X}\to \mathcal{X}$ 
   on a compact metric space $(\mathcal{X}, \mathbf{d})$ with invariant probability measures 
   $\nu_1, \nu_2 \in \mathscr{M}^T(\mathcal{X})$ satisfying 
   \begin{itemize}
     \item the upper metric mean dimension $\umdimm\left(\mathcal{X}, T, \mathbf{d}\right)$ is finite,
     \item $\urdim\left(\mathcal{X}, T, \mathbf{d}, \nu_1 \right) = \urdim\left(\mathcal{X}, T, \mathbf{d}, \nu_2 \right) =1$,
     \item $\urdim\left(\mathcal{X}, T, \mathbf{d}, \frac{1}{2}\nu_1 + \frac{1}{2}\nu_2 \right) \leq \frac{1}{2}$.
   \end{itemize}
\end{proposition}
This shows that the inequality (\ref{eq: convexity inequality for upper rate distortion dimension}) may be strict in general.
(Consider (\ref{eq: convexity inequality for upper rate distortion dimension}) 
for $\lambda := \frac{1}{2}\delta_{\nu_1} + \frac{1}{2}\delta_{\nu_2}$.)

Proposition \ref{proposition: convexity inequality may be strict} follows from 

\begin{proposition} \label{proposition: exotic example}
 Let $a_1<b_1<a_2<b_2<a_3<b_3<\dots$ be an increasing sequence of natural numbers with $b_k > k^2 a_k$ for all $k\geq 1$.
 There exists a continuous action $T\colon \mathbb{R}^d\times \mathcal{X}\to \mathcal{X}$ 
 on a compact metric space $(\mathcal{X}, \mathbf{d})$ 
 with an invariant probability measure $\nu \in \mathscr{M}^T(\mathcal{X})$ satisfying 
 \begin{itemize}
   \item  the upper metric mean dimension $\umdimm\left(\mathcal{X}, T, \mathbf{d}\right)$ is finite,
   \item  $\urdim \left(\mathcal{X}, T, \mathbf{d}, \nu\right) = 1$,
   \item  ${\displaystyle \lim_{\substack{t\in \bigcup_{k} [b_k, a_{k+1}) \\ t\to \infty}}\frac{R(\mathbf{d}, \nu, 2^{-t})}{t} = 0}$, where 
   the real parameter $t$ goes to infinity while it satisfies $t\in \bigcup_{k=1}^\infty [b_k, a_{k+1})$.
 \end{itemize}
\end{proposition}

We first prove Proposition \ref{proposition: convexity inequality may be strict}, assuming Proposition \ref{proposition: exotic example}.

\begin{proof}[Proof of Proposition \ref{proposition: convexity inequality may be strict}]
We take an increasing sequence of natural numbers $a_1<b_1<a_2<b_2<\dots$ with 
$b_k>k^2 a_k$ and $a_{k+1}>k^2 b_k$.
We set 
\[ a^{(1)}_k = a_k,\quad b^{(1)}_k = b_k,\quad 
    a^{(2)}_k = b_k, \quad b^{(2)}_k = a_{k+1} \quad (k\geq 1). \]
We have $b^{(i)}_k > k^2 a^{(i)}_k$ for both $i=1, 2$.
For $i=1, 2$,  by applying Proposition \ref{proposition: exotic example} to the sequence 
$a^{(i)}_1<b^{(i)}_1< a^{(i)}_2 < b^{(i)}_2 < \dots$, we find a
continuous action $T_i \colon \mathbb{R}^d\times \mathcal{X}_i \to \mathcal{X}_i$ 
 on a compact metric space $(\mathcal{X}_i, \mathbf{d}^{(i)})$ 
 with an invariant probability measure $\nu_i \in \mathscr{M}^{T_i}(\mathcal{X}_i)$ satisfying 
\begin{itemize}
   \item the upper metric mean dimension $\umdimm\left(\mathcal{X}_i, T_i, \mathbf{d}^{(i)} \right)$ is finite,
   \item $\urdim \left(\mathcal{X}_i, T_i, \mathbf{d}^{(i)}, \nu_i \right) = 1$,
   \item  ${\displaystyle \lim_{\substack{t\in \bigcup_{k} [b^{(i)}_k, a^{(i)}_{k+1}) \\ t\to \infty}}
             \frac{R(\mathbf{d}^{(i)}, \nu_i, 2^{-t})}{t} = 0}$.
\end{itemize}

We set $\mathcal{X} = \mathcal{X}_1\cup \mathcal{X}_2$ (the disjoint union).
We define a metric $\mathbf{d}$ on it by 
\[ \mathbf{d}(x, y) = \begin{cases} \mathbf{d}^{(1)}(x, y) & (x, y\in \mathcal{X}_1) \\
                                                 \mathbf{d}^{(2)}(x, y) & (x, y\in \mathcal{X}_2) \\
                                                 \max\left(\diam(\mathcal{X}_1, \mathbf{d}^{(1)}), \diam(\mathcal{X}_2, \mathbf{d}^{(2)})\right)
                                                 & (\text{otherwise})  \end{cases}. \]
We define $T\colon \mathbb{R}^d\times \mathcal{X}\to \mathcal{X}$ as 
$T^u x = T_i^u x$ for $x\in \mathcal{X}^{(i)}$ and $u\in \mathbb{R}^d$.
We can naturally think $\nu_i$ as invariant probability measures on $\mathcal{X}$.
It is immediate to see that 
  \begin{itemize}
    \item the upper metric mean dimension $\umdimm\left(\mathcal{X}, T, \mathbf{d}\right)$ is finite, 
    \item $\rdim\left(\mathcal{X}, T, \mathbf{d}, \nu_i \right) = \rdim\left(\mathcal{X}_i, T_i, \mathbf{d}^{(i)}, \nu_i\right) = 1$.
  \end{itemize}
Therefore we only need to check that $\urdim\left(\mathcal{X}, T, \mathbf{d}, \frac{1}{2}\nu_1 + \frac{1}{2}\nu_2\right) \leq \frac{1}{2}$.
Set $\mu = \frac{1}{2}\nu_1 + \frac{1}{2}\nu_2$.

Let $t$ be a large real number.
By Theorem \ref{theorem: rate distortion function and convex combination of measures}
\[  R(\mathbf{d}, \mu, 2^{-t}) \leq \frac{1}{2} R(\mathbf{d}^{(1)}, \nu_1, 2^{-t}) + \frac{1}{2} R(\mathbf{d}^{(2)}, \nu_2, 2^{-t}). \]
Here $R(\mathbf{d}^{(i)}, \nu_i, 2^{-t})$ denotes the rate distortion function of $(\mathcal{X}_i, T_i, \mathbf{d}^{(i)}, \nu_i)$.
It follows from $\rdim\left(\mathcal{X}_i, T_i, \mathbf{d}^{(i)}, \nu_i\right) = 1$ 
that $R(\mathbf{d}^{(i)}, \nu_i, 2^{-t}) = t + o(t)$ for both $i=1, 2$, where the error term satisfies $\lim_{t\to \infty} o(t)/t = 0$.

If $t \in [b_k, a_{k+1}) = [b^{(1)}_k, a^{(1)}_{k+1})$ for some $k\geq 1$, then $R(\mathbf{d}^{(1)}, \nu_1, 2^{-t}) = o(t)$.
If $t \in [a_{k+1}, b_{k+1}) = [b^{(2)}_k, a^{(2)}_{k+1})$ for some $k\geq 1$, then 
$R(\mathbf{d}^{(2)}, \nu_2, 2^{-t}) = o(t)$.
Therefore, in both cases, we have 
\[  R(\mathbf{d}, \mu, 2^{-t})  \leq  \frac{1}{2}(t+o(t)) + \frac{1}{2}o(t) = \frac{t}{2} + o(t). \]
Thus we conclude 
\[  \urdim\left(\mathcal{X}, T, \mathbf{d}, \mu\right) = \limsup_{t\to \infty} \frac{R(\mathbf{d}, \mu, 2^{-t})}{t} \leq \frac{1}{2}. \]
\end{proof}

The rest of this section is devoted to the proof of Proposition \ref{proposition: exotic example}.
Our construction is similar to that of \S \ref{section: construction of an example}, and we
use the terminology “$\Gamma$-piecewisely linear map” introduced there.

\begin{figure}[h] 
    \centering
    \includegraphics[width=3.0in]{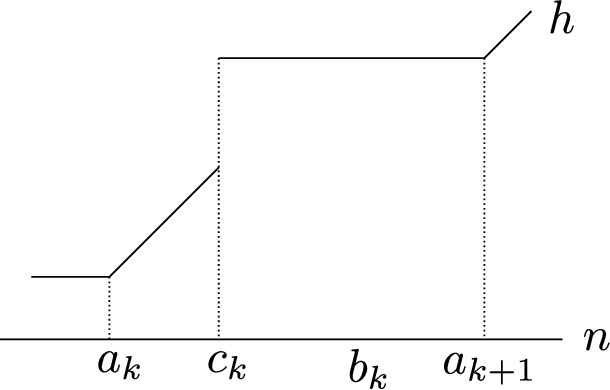}
    \caption{The graph of the function $h$. We have $c_k = k a_k$ and $k c_k < b_k < a_{k+1}$.}    \label{figure: function_h}
\end{figure}

Suppose that we are given an increasing sequence of natural numbers $a_1<b_1<a_2<b_2<a_3<b_3<\dots$ satisfying 
$b_k > k^2 a_k$.
Set $c_k = k a_k$. We have $a_k \leq c_k \leq b_k$.
We define a function $h\colon \mathbb{N}\to \mathbb{N}$ by 
\[  h(n) = \begin{cases}  n & (n\leq  a_1 \text{ or } n\in \bigcup_{k=1}^\infty [a_k, c_k]) \\
                                  a_{k+1} & (n\in  (c_k, a_{k+1}] \text{ for $k\geq 1$}) \end{cases}. \]
Figure \ref{figure: function_h} shows a schematic picture of the graph of the function $h$. 
Notice that $h$ is monotone non-decreasing and satisfies $h(n) \geq n$ for all $n$.

Let $(V, \norm{\cdot})  = (\ell^\infty(\mathbb{N}), \norm{\cdot}_{\ell^\infty})$ be the Banach space of 
bounded sequences (indexed by natural numbers).
Let $B = \{v\in V\mid \norm{v}\leq 1\}$ be the unit ball.
We define $C(\mathbb{R}^d, B)$ as the space of all continuous maps $x\colon \mathbb{R}^d\to B$.
We define a metric $\mathbf{d}$ on it by 
\[ \mathbf{d}(x, x^\prime) = \sup_{n\geq 1} \left(2^{-n}\sup_{|u|_\infty \leq n} \norm{x(u)-x^\prime(u)}\right). \]
Let $T\colon \mathbb{R}^d\times  C(\mathbb{R}^d, B)\to C(\mathbb{R}^d, B)$ be the natural shift action ($T^u x(t) = x(t+u)$).

Let $\Lambda = \{0, 1\}^{\mathbb{N}}$ be the full shift on the alphabet $0,1$.
We define a metric $\rho$ on it by
\[ \rho(v, w)  =  2^{-h\left(\min\{m\geq 1\mid  v_m\neq w_m\}\right)}. \] 
Since $h(n)\geq n$, we have $\rho(v, w) \leq 2^{-\min\{m\geq 1\mid v_m\neq w_m\}}$ where the right-hand side is a more 
standard metric on $\Lambda$.
We define a map $\varphi\colon \Lambda \to B$ by 
\[  \varphi(v) = \left(2^{-h(m)}v_m\right)_{m=1}^\infty. \]
Set $\Lambda^\prime = \varphi(\Lambda)$. 
The metric spaces $(\Lambda, \rho)$ and $(\Lambda^\prime, \norm{\cdot})$ are isometric, namely 
\[ \norm{\varphi(v)-\varphi(w)} = \rho(v, w). \]
In particular, $\diam\left(\Lambda^\prime, \norm{\cdot}\right) = \diam\left(\Lambda, \rho\right)  = \frac{1}{2}$.

As in \S \ref{section: construction of an example} we fix a $\mathbb{Z}^d$-invariant triangulation $\Gamma$ of $\mathbb{R}^d$ such that 
its vertex set is equal to $\mathbb{Z}^d$.
We define $\mathcal{Y}\subset C(\mathbb{R}^d, B)$ as the space of all $\Gamma$-piecewisely linear maps 
$x\colon \mathbb{R}^d\to B$ satisfying $x(\mathbf{n}) \in \Lambda^\prime$ for all $\mathbf{n}\in \mathbb{Z}^d$.
We define $\mathcal{X}\subset C(\mathbb{R}^d, B)$ by 
\[ \mathcal{X} = \bigcup_{s\in [0, 1)^d} T^s \mathcal{Y}. \]
This is a compact $T$-invariant subset of $C(\mathbb{R}^d, B)$.

Let $p = \left(\frac{1}{2}\delta_0 + \frac{1}{2}\delta_1\right)^{\otimes \mathbb{N}}$ be the (unbiased) Bernoulli measure on 
$\Lambda = \{0, 1\}^{\mathbb{N}}$.
Let $p^{\otimes \mathbb{Z}^d}$ be the product measure on $\Lambda^{\mathbb{Z}^d}$, and
let $\Phi\colon \Lambda^{\mathbb{Z}^d}\to \mathcal{Y}$ be the natural map.
We set $\mu = \Phi_* \left(p^{\otimes \mathbb{Z}^d}\right)$. 
This is a probability measure on $\mathcal{Y}$.
We define $\nu \in \mathscr{M}^T(\mathcal{X})$ by 
\[ \nu = \int_{[0, 1)^d} T^s_*\mu\, d\mathbf{m}(s). \]
We will prove that $(\mathcal{X}, T, \mathbf{d}, \nu)$ satisfies the claim of Proposition \ref{proposition: exotic example}.

\begin{lemma} \label{lemma: upper metric mean dimension and gaped limit}
  \begin{enumerate}
    \item  $\umdimm\left(\mathcal{X}, T, \mathbf{d}\right) \leq 1$. 
              In particular, $\urdim\left(\mathcal{X}, T, \mathbf{d}, \nu\right) \leq 1$.
    \item  ${\displaystyle \lim_{\substack{t\in \bigcup_{k} [b_k, a_{k+1}) \\ t\to \infty}}\frac{R(\mathbf{d}, \nu, 2^{-t})}{t} = 0}$.
  \end{enumerate}
\end{lemma}

\begin{proof}
As in Lemma \ref{lemma: Lipschitz constant of PL maps}, we can find $C>0$ such that every $x\in \mathcal{X}$
is $C$-Lipschitz. 
In particular, for any $L>0$ and $s, t\in \mathbb{R}^d$
\[ \mathbf{d}_L(T^s x, T^t x) \leq \frac{C}{2} |s-t|_\infty. \]

(1) For $L>0$ and $\varepsilon>0$
\[ \#\left(\mathcal{X}, \mathbf{d}_L, \varepsilon\right) \leq 
   \sum_{s\in [0, 1)^d\cap \left(\frac{\varepsilon}{2C}\mathbb{Z}^d\right)} \#\left(T^s\mathcal{Y}, \mathbf{d}_L, \frac{\varepsilon}{2}\right) 
   \leq \left(1+\frac{2C}{\varepsilon}\right)^d \#\left(\mathcal{Y}, \mathbf{d}_{L+1}, \frac{\varepsilon}{2}\right). \]
Choose $\ell = \ell(\varepsilon)>0$ satisfying $2^{-\ell} < \frac{\varepsilon}{2}$.
Then
\[  \#\left(\mathcal{Y}, \mathbf{d}_{L+1}, \frac{\varepsilon}{2}\right) 
     \leq  \#\left(\Lambda, \rho, \varepsilon\right)^{(L+2\ell+1)^d}. \]
Hence 
\[ S(\mathcal{X}, T, \mathbf{d}, \varepsilon) = \lim_{L\to \infty} \frac{\log  \#\left(\mathcal{X}, \mathbf{d}_L, \varepsilon\right)}{L^d}
    \leq \log \#\left(\Lambda, \rho, \varepsilon\right). \]
Therefore 
\[ \umdimm\left(\mathcal{X}, T, \mathbf{d}\right) 
  = \limsup_{\varepsilon\to 0} \frac{S(\mathcal{X}, T, \mathbf{d}, \varepsilon)}{\log (1/\varepsilon)}
  \leq  \limsup_{\varepsilon\to 0} \frac{\log \#\left(\Lambda, \rho, \varepsilon\right)}{\log (1/\varepsilon)}. \]
The right-most side is the upper Minkowski dimension of $(\Lambda, \rho)$.
Noting $h(n)\geq n$, it is easy to see that the upper Minkowski dimension is less than or equal to one.

(2) We again use the inequality $S(\mathcal{X}, T, \mathbf{d}, \varepsilon) \leq \log \#\left(\Lambda, \rho, \varepsilon\right)$.
Let $\varepsilon = 2^{-t}$ with $b_k\leq t < a_{k+1}$. 
Consider the natural projection $\Lambda\to \{0, 1\}^{c_k}$ to the first $c_k$ coordinates.
It follows from the definition of the metric $\rho$ that 
its fibers have diameter $2^{-a_{k+1}} < \varepsilon$. 
Hence we have 
\[ \log \#\left(\Lambda, \rho, \varepsilon\right) \leq c_k < \frac{b_k}{k} \leq \frac{t}{k}. \]
Namely 
\[ \frac{R(\mathbf{d}, \nu, 2^{-t})}{t} \leq  \frac{S(\mathcal{X}, T, \mathbf{d}, 2^{-t})}{t}  \leq \frac{1}{k} \quad (b_k\leq t < a_{k+1}). \]
We have $k\to \infty$ as $t\to \infty$. So this shows the claim.
\end{proof}

Now we only need to show that $\urdim(\mathcal{X}, T, \mathbf{d}, \nu) = 1$.
In fact, we will prove that $\lim_{k\to \infty} \frac{R(\mathbf{d}, \nu, 2^{-{c_k}})}{c_k} = 1$.
Then $\urdim(\mathcal{X}, T, \mathbf{d}, \nu) = 1$ follows from Lemma \ref{lemma: upper metric mean dimension and gaped limit} (1).

As in \S \ref{section: construction of an example}, we would like to use Proposition \ref{proposition: lower bound on rate distortion function}.
For that purpose we need to estimate 
\[   \int_{\mathcal{X}} 2^{-\frac{a}{L^d}\int_{[0, L)^d} \mathbf{d}(T^t x, y_t)d\mathbf{m}(t)} d\nu(x)
     = \int_{[0,1)^d}\left(\int_{T^s \mathcal{Y}}2^{-\frac{a}{L^d}\int_{[0, L)^d} \mathbf{d}(T^t x, y_t)d\mathbf{m}(t)} dT^s_*\mu(x)\right) d\mathbf{m}(s) \]
for $L>0$, $y\in L^1([0,L)^d, \mathcal{X})$ and appropriately chosen $a>0$.
(The choice $a=2^{c_k+2}L^d$ will work.)

Let $L\in \mathbb{N}$, $s\in [0,1)^d$ and $y\in L^1([0, L)^d, \mathcal{X})$.
By the same argument as in Lemma \ref{lemma: distortion estimate}, for every $x\in \mathcal{X}$ we have 
\[ \int_{[0, L)^d}\mathbf{d}(T^t x, y_t)\, d\mathbf{m}(t) \geq \frac{1}{2}\sum_{\mathbf{n}\in \mathbb{Z}^d\cap [0, L)^d}
    \norm{x(s+\mathbf{n})-\int_{\mathbf{n}+[0,1)^d} y_t(s+\mathbf{n}-t)\, d\mathbf{m}(t)}. \]
Let 
\[ z_s(\mathbf{n}) = \int_{\mathbf{n}+[0, 1)^d} y_t(s+\mathbf{n}-t)\, d\mathbf{m}(t) \in B  \quad 
    (\mathbf{n}\in \mathbb{Z}^d\cap [0, L)^d). \]
For any $a>0$
\begin{equation} \label{eq: estimate of partition-like function}
  \begin{split}
  \int_{T^s \mathcal{Y}}2^{-\frac{a}{L^d}\int_{[0, L)^d} \mathbf{d}(T^t x, y_t)d\mathbf{m}(t)} dT^s_*\mu(x) & \leq 
  \int_{T^s\mathcal{Y}} 2^{-\frac{a}{2L^d}\sum_{\mathbf{n}\in \mathbb{Z}^d\cap [0, L)^d} \norm{x(s+\mathbf{n})-z_s(\mathbf{n})}}dT^s_*\mu(x) \\
  & =  \int_{T^s\mathcal{Y}} \prod_{\mathbf{n}\in \mathbb{Z}^d\cap [0, L)^d} 
         2^{-\frac{a}{2L^d}\norm{x(s+\mathbf{n})-z_s(\mathbf{n})}} dT^s_*\mu(x) \\
  & = \prod_{\mathbf{n}\in \mathbb{Z}^d\cap [0, L)^d}  \int_{\Lambda} 2^{-\frac{a}{2L^d}\norm{\varphi(v)-z_s(\mathbf{n})}} dp(v).     
  \end{split}
\end{equation}
We would like to use this for $a=2^{c_k+2}L^d$.

\begin{lemma} \label{lemma: estimate of integral over full-shift}
There exists a universal positive constant $K$ such that for any $z\in V$ and 
sufficiently large $k$ (independent of $z$) we have 
\[ \int_{\Lambda} 2^{-2^{c_k+1}\norm{\varphi(v)-z}} dp(v) \leq  K 2^{-c_k}. \]
\end{lemma}

\begin{proof}
Take $z\in V$.
Consider a function $\Lambda\ni v\mapsto \norm{\varphi(v)-z}\in \mathbb{R}$, and let $w = (w_n)_{n=1}^\infty \in \Lambda$ be a point 
attaining the minimum of this function.
Then for any point $v\in \Lambda$ we have
\[ \rho(v, w) = \norm{\varphi(v)-\varphi(w)} \leq \norm{\varphi(v)-z} + \norm{z-\varphi(w)} \leq 2\norm{\varphi(v)-z}. \]
Hence $\int_{\Lambda} 2^{-2^{c_k+1}\norm{\varphi(v)-z}} dp(v)  \leq  \int_{\Lambda} 2^{-2^{c_k}\rho(v, w)} dp(v)$.

We define subsets $\Lambda_i$ ($i=1, \dots, c_k$) of $\Lambda$ by 
$\Lambda_1 = \{v \mid v_1\neq w_1\}, \Lambda_2 = \{v\mid v_1=w_1, v_2\neq w_2\}, \Lambda_3 = \{v\mid v_1=w_1, v_2 = w_2, v_3\neq w_2\},
\dots, \Lambda_{c_k} = \{v\mid v_1=w_1, \dots, v_{c_k-1}=w_{c_k-1}, v_{c_k}\neq w_{c_k}\}$.
We also set 
\[ \Lambda^\prime = \Lambda\setminus (\Lambda_1\cup \dots, \cup \Lambda_{c_k}) 
     = \{v\mid v_1=w_1, \dots, v_{c_k}=w_{c_k}\}. \]
Then 
\begin{align*}
  & \int_{\Lambda} 2^{-2^{c_k}\rho(v, w)} dp(v) \\ & = \int_{\bigcup_{i=1}^{a_k}\Lambda_i} 2^{-2^{c_k}\rho(v, w)} dp(v)
   + \sum_{i=a_k+1}^{c_k} \int_{\Lambda_i} 2^{-2^{c_k}\rho(v, w)} dp(v) + \int_{\Lambda^\prime} 2^{-2^{c_k}\rho(v, w)} dp(v). 
\end{align*}   
We have $\rho(v, w) \geq 2^{-a_k}$ over $\bigcup_{i=1}^{a_k}\Lambda_i$.
On $\Lambda_i$ $(a_k+1\leq i \leq c_k)$ we have $\rho(v, w) = 2^{-i}$.
We also have $p(\Lambda^\prime) = 2^{-c_k}$. Hence 
\[ \int_{\Lambda} 2^{-2^{c_k}\rho(v, w)} dp(v) \leq  2^{-2^{c_k-a_k}} + \sum_{i=a_k+1}^{c_k} 2^{-i} 2^{-2^{c_k-i}} + 2^{-c_k}. \]
Since $c_k = k a_k$, we have $2^{-2^{c_k-a_k}} < 2^{-c_k}$ for large $k$.
For estimating the second term, we introduce a change of variable $j = c_k-i$:
\[ \sum_{i=a_k+1}^{c_k} 2^{-i} 2^{-2^{c_k-i}} = \sum_{j=0}^{c_k-a_k-1}2^{-c_k+j} 2^{-2^j} 
    < 2^{-c_k} \sum_{j=0}^\infty 2^j 2^{-2^j}. \]
Then
\[ \int_{\Lambda} 2^{-2^{c_k}\rho(v, w)} dp(v) \leq  \left(2+ \sum_{j=0}^\infty 2^j 2^{-2^j}\right) 2^{-c_k}. \]
Now the claim holds for $K := 2+ \sum_{j=0}^\infty 2^j 2^{-2^j} < \infty$.
\end{proof}

We assume $k$ is large.
By applying Lemma \ref{lemma: estimate of integral over full-shift} to (\ref{eq: estimate of partition-like function})
with $a=2^{c_k+2}L^d$ 
\[  \int_{T^s \mathcal{Y}} 2^{-2^{c_k+2}\int_{[0, L)^d} \mathbf{d}(T^t x, y_t)d\mathbf{m}(t)} dT^s_*\mu(x)
    \leq K^{L^d} 2^{-c_k L^d}. \]
By integrating this over $s\in [0, 1)^d$,
\[ \int_{\mathcal{X}} 2^{-2^{c_k+2}\int_{[0, L)^d} \mathbf{d}(T^t x, y_t)d\mathbf{m}(t)} d\nu(x) 
    \leq  K^{L^d} 2^{-c_k L^d}. \]
Now we use Proposition \ref{proposition: lower bound on rate distortion function}
with the parameters $\varepsilon = 2^{-c_k}$, $a=2^{c_k+2} L^d$ and $\lambda = K^{-L^d} 2^{c_k L^d}$.
Then
\[ R_L(\mathbf{d}, \nu, 2^{-c_k}) \geq -a\varepsilon + \log \lambda  = -4L^d - L^d \log K + c_k L^d. \]
Dividing this by $L^d$ and letting $L\to \infty$, we get 
\[  R(\mathbf{d}, \nu, 2^{-c_k}) \geq  -4 - \log K + c_k. \]

\begin{lemma} \label{lemma: rate distortion dimension along sparse parameters} 
$\lim_{k\to \infty} \frac{R(\mathbf{d}, \nu, 2^{-c_k})}{c_k} = \urdim(\mathcal{X}, T, \mathbf{d}, \nu) =1$.
\end{lemma}

\begin{proof}
Assume $k$ is large.
From $R(\mathbf{d}, \nu, 2^{-c_k}) \geq  -4 - \log K + c_k$ 
\[ 1 \leq \liminf_{k\to \infty} \frac{R(\mathbf{d}, \nu, 2^{-c_k})}{c_k} \leq
    \limsup_{k\to \infty} \frac{R(\mathbf{d}, \nu, 2^{-c_k})}{c_k} \leq \urdim(\mathcal{X}, T, \mathbf{d}, \nu) \leq 1. \]
The last inequality is given by Lemma \ref{lemma: upper metric mean dimension and gaped limit} (1).
\end{proof}

By combining Lemmas \ref{lemma: upper metric mean dimension and gaped limit} 
and \ref{lemma: rate distortion dimension along sparse parameters}, we have proved Proposition \ref{proposition: exotic example}.

\end{document}